\documentclass[11pt, a4paper]{amsart}
\usepackage{mathbbol}
\usepackage{amsmath, amsthm, amssymb}
\usepackage{txfonts, mathrsfs}
\usepackage[usenames]{color}
\usepackage{psfrag}
\usepackage{graphicx}

\newtheorem*{thm*}{Theorem}

\newtheorem{thm}{Theorem}[section]
\newcommand{\bt}{\begin{thm}}
\newcommand{\et}{\end{thm}}

\newtheorem{cor}[thm]{Corollary}
\newcommand{\bc}{\begin{cor}}
\newcommand{\ec}{\end{cor}}

\newtheorem{lem}[thm]{Lemma}
\newcommand{\bl}{\begin{lem}}
\newcommand{\el}{\end{lem}}

\newtheorem{prop}[thm]{Proposition}
\newcommand{\bp}{\begin{prop}}
\newcommand{\ep}{\end{prop}}

\newtheorem{defn}[thm]{Definition}
\newcommand{\bd}{\begin{defn}}      
\newcommand{\ed}{\end{defn}}

\newtheorem{rmrk}[thm]{Remark}
\newcommand{\br}{\begin{rmrk}}
\newcommand{\er}{\end{rmrk}}

\newtheorem{quest}[thm]{Question}
\newcommand{\bq}{\begin{quest}}
\newcommand{\eq}{\end{quest}}

\newtheorem{example}[thm]{Example}

\newcommand{\N}{\mathbb{N}}

\newcommand{\R}{\mathbb{R}}

\newdimen\vintkern\vintkern12pt
\def\vint{-\kern-\vintkern\int}

\newcommand{\hm}{{\mathcal H}}

\newcommand{\dist}{\operatorname{dist}}
\newcommand{\diam}{\operatorname{diam}}
\newcommand{\trace}{\operatorname{tr}}
\newcommand{\length}{\ell}
\newcommand{\Area}{\operatorname{Area}}

\newcommand{\md}{\operatorname{md}}

\newcommand{\jac}{{\mathbf J}}
\newcommand{\ap}{\operatorname{ap}}
\newcommand{\apmd}{\ap\md}

\DeclareMathOperator{\MOD}{mod}

\begin{document}
\bibliographystyle{plain}

\title[]{Canonical parametrizations of metric discs}

\author{Alexander Lytchak}


\address
  {Mathematisches Institut\\ Universit\"at K\"oln\\ Weyertal 86 -- 90\\ 50931 K\"oln, Germany}
\email{alytchak@math.uni-koeln.de}

\author{Stefan Wenger}

\address
  {Department of Mathematics\\ University of Fribourg\\ Chemin du Mus\'ee 23\\ 1700 Fribourg, Switzerland}
\email{stefan.wenger@unifr.ch}

\date{\today}

\thanks{S.~W.~was partially supported by Swiss National Science Foundation Grants 153599 and 165848.}

\begin{abstract}
We use the recently established existence and regularity of area and energy minimizing discs in metric spaces to obtain canonical parametrizations of metric surfaces. Our approach yields a new and conceptually simple proof of a well-known theorem of Bonk and Kleiner on the existence of quasisymmetric parametrizations of linearly locally connected, Ahlfors $2$-regular metric $2$-spheres. Generalizations and applications to the geometry of such surfaces are described.
\end{abstract}

\maketitle

\section{Introduction and statement of main results}

By the classical uniformization theorem, every smooth Riemann surface is conformally diffeomorphic to a surface of constant curvature. A fundamental question, widely studied in the field of Analysis in Metric Spaces, asks to what extent non-smooth metric spaces admit parametrizations from a Euclidean domain with good geometric and analytic properties.  In this more general context, one usually looks for parametrizations which are biLipschitz, quasisymmetric, or quasiconformal.

 A celebrated and difficult theorem of Bonk--Kleiner \cite{BK02} asserts that an Ahlfors $2$-regular metric space $X$, homeomorphic to the standard $2$-sphere $S^2$, admits a quasisymmetric homeomorphism from $S^2$ to $X$ if and only if $X$ is linearly locally connected. We refer to Section~\ref{sec:prelims} for the definitions of quasisymmetric homeomorphism, linear local connectedness, and Ahlfors $2$-regularity. This result has since been extended for example in \cite{BK05}, \cite{Wil08}, \cite{Wil10}, \cite{MW13}, \cite{LRR17} and recently in the important paper \cite{Raj14}. We refer to \cite{Raj14} for details and more references.

The main purpose of the present paper is to provide a new and conceptually simple approach to the  theorem of Bonk--Kleiner and related results. Our approach  is a direct generalization  of the classical existence proof of conformal  parametrizations of smooth surfaces via minimizing the energy of maps into the surface, see \cite{Jos91-book}. Our main parametrization result can be stated as follows. We denote by $D$ and $\overline{D}$ the open and closed unit discs in $\R^2$, respectively, and refer to Section~\ref{sec:prelims} for the definition of energy $E_+^2(u)$ of a (Sobolev) map $u$ from $D$ to a metric space.

\bt\label{thm:main-thm-qs-discs-summary}
 Let $X$ be a geodesic metric space homeomorphic to $\overline{D}$ and with boundary circle $\partial X$ of finite length. If $X$ is Ahlfors $2$-regular and linearly locally connected then 
there exists a homeomorphism $u\colon\overline{D}\to X$ of minimal energy $E_+^2(u)<\infty$.
Any such $u$ is quasisymmetric and is uniquely determined up to a conformal diffeomorphism of $\overline{D}$.
\et

A more general statement will be provided in Theorem~\ref{thm:qs-subdomains-llc}. 
Note that the theorem comprises several statements, which will be described more precisely and in a different order below. Firstly, there exists a continuous map of finite energy from $\overline{D}$ to $X$ whose boundary parametrizes the boundary circle $\partial X$. Secondly, there exists one such map of minimal energy. Finally, any such map is a quasisymmetric homeomorphism which is unique up to composition with a conformal diffeomorphism of $\overline{D}$.
Similarly to Theorem~\ref{thm:main-thm-qs-discs-summary}, we obtain a canonical (up to conformal diffeomorphism) quasisymmetric parametrization of Ahlfors $2$-regular, linearly locally connected metric $2$-spheres, see Theorem~\ref{thm:BK-canonical-version}, and consequently the Bonk--Kleiner theorem mentioned above, see Corollary~\ref{cor:Bonk-Kleiner} below.

We now describe our approach and the statements of the theorem in more precise terms. Let $X$ be a complete metric space. Denote by $N^{1,2}(D, X)$ the space of (Newton-) Sobolev maps from $D$ to $X$ in the sense of \cite{HKST15}. Given a map $u\in N^{1,2}(D, X)$ we let $\trace(u)$ be its trace, $E_+^2(u)$ its Reshetnyak energy, and $\Area(u)$ its parametrized area. See Section~\ref{sec:Sobolev-defs} for these definitions and references. If $\Gamma\subset X$ is a Jordan curve then we let $\Lambda(\Gamma, X)$ be the possibly empty family of maps $u\in N^{1,2}(D, X)$ whose trace has a continuous representative which is a weakly monotone parametrization of $\Gamma$. 

Our first result provides topological information on energy minimizers:

\bt\label{thm:Plateau-solution-cell-like-intro}
 Let $X$ be a geodesic metric space homeomorphic to $\overline{D}$ and let $u\colon\overline{D}\to X$ be a continuous map. If  $u$ is in $\Lambda(\partial X, X)$ and minimizes the Reshetnyak energy $E_+^2$ among all maps in $\Lambda(\partial X,X)$ then $u$ is a uniform limit of homeomorphisms from $\overline{D}$ to $X$.
\et

 If $X$ has a quadratic bound $\hm^2(B(x,r))\leq C\cdot r^2$ for the Hausdorff $2$-measure of $r$-balls then any $u$ as in the theorem above is a homeomorphism, see Theorem~\ref{thm:monotone-implies-homeo-quadgrowth}. 

In general, the family $\Lambda(\partial X, X)$ may be empty for spaces as in Theorem~\ref{thm:Plateau-solution-cell-like-intro}. However, for spaces admitting a quadratic isoperimetric inequality as defined below, this family is not empty when the length of $\partial X$ is finite. Furthermore, energy minimizers exist, are continuous, and  their infinitesimal structure is as close to conformal as possible by our results in \cite{LW15-Plateau}, \cite{LW16-energy-area}, see also Theorem~\ref{thm:existence-regularity-energy-min} below.

\bd\label{def:quad-isop}
 A complete metric space $X$ is said to admit a quadratic isoperimetric inequality if there exists $C>0$ such that every Lipschitz curve $c\colon S^1\to X$ is the trace of some $u\in N^{1,2}(D, X)$ such that $$\Area(u)\leq C\cdot \length(c)^2,$$ where $\length(c)$ denotes the length of $c$.
\ed

The isoperimetric constant of $X$ is defined as the infimum over all $C>0$ for which the above holds. 
Our second ingredient in the proof of Theorem~\ref{thm:main-thm-qs-discs-summary} is:

\bt\label{thm:isop-ineq-for-geod-discs}
 Let $X$ be a complete, geodesic metric space homeomorphic to a $2$-dimensional manifold. Suppose there exists $C>0$ such that every Jordan curve in $X$ bounds a Jordan domain $U\subset X$ with 
\begin{equation}\label{eq:geom-isop-intro}
\hm^2(U)\leq C\cdot \length(\partial U)^2.
\end{equation}
Then $X$ admits a quadratic isoperimetric inequality. Moreover, the isoperimetric constant of $X$ is at most $C$.  
\et

By a manifold we mean a topological manifold with or without boundary. A Jordan domain $U\subset X$ is an open set homeomorphic to $D$ such that $\overline{U}\subset X$ is homeomorphic to $\overline{D}$.
 %In particular, spaces $X$ as in the theorem above are necessarily simply connected. 
Particular examples of spaces satisfying \eqref{eq:geom-isop-intro} are Ahlfors $2$-regular, linearly locally connected metric spaces homeomorphic to $\overline{D}$ or $S^2$, see Lemma~\ref{lem:linloccon-settheoretic-isop}.

%
%In \cite{LW15-Plateau}, \cite{LW16-energy-area}, we proved existence, regularity, and an infinitesimal quasiconformality property for energy minimizers in $\Lambda(\Gamma, X)$ whenever $\Gamma$ has finite length and $X$ is proper and admits a quadratic isoperimetric inequality. See the summarizing Theorem~\ref{thm:existence-regularity-energy-min} below. In view of Theorem~\ref{thm:isop-ineq-for-geod-discs}, these results apply to spaces as in Theorem~\ref{thm:main-thm-qs-discs-summary}.
 
Now, if $X$ is as in Theorem~\ref{thm:main-thm-qs-discs-summary} then it admits a quadratic isoperimetric inequality by Theorem~\ref{thm:isop-ineq-for-geod-discs}. Hence, the family $\Lambda(\partial X, X)$ is not empty and, by \cite{LW15-Plateau}, \cite{LW16-energy-area}, contains an energy minimizer $u$ which is continuous. By Theorem~\ref{thm:Plateau-solution-cell-like-intro} and the remark following it, any such $u$ is a homeomorphism. 
Since the infinitesimal structure of energy minimizers is as close to conformal as possible by \cite{LW15-Plateau}, \cite{LW16-energy-area}, the modulus estimates from \cite{HK98} imply that $u$ is
quasisymmetric; furthermore $u$ is uniquely determined up to a conformal diffeomorphism of $\overline{D}$. This concludes the outline of the proof of 
Theorem~\ref{thm:main-thm-qs-discs-summary}.

The class of metric surfaces satisfying \eqref{eq:geom-isop-intro} is of importance in the theory of minimal surfaces in metric spaces, due to \cite{LW16-intrinsic}. Note that such surfaces need not be Ahlfors $2$-regular and they need not admit quasiconformal parametrizations (as defined in \cite{Raj14}) even locally, see Example~\ref{ex:collapsed-disc}. Combining Theorems~\ref{thm:isop-ineq-for-geod-discs} and \ref{thm:Plateau-solution-cell-like-intro} we nevertheless obtain natural "almost parametrizations" of such surfaces, see Theorem~\ref{thm:almost-qc-param}.
For surfaces satisfying \eqref{eq:geom-isop-intro} with the Euclidean constant $C=\frac{1}{4\pi}$ we obtain, as a consequence of Theorem~\ref{thm:isop-ineq-for-geod-discs} together with \cite[Theorem 1.1]{LW-isoperimetric}, the following strengthening of \cite[Theorem 1.3]{LW-isoperimetric}. 

\bc
 Let $X$ be a geodesic metric space homeomorphic to $\overline{D}$. Then $X$ is a ${\rm CAT}(0)$-space if and only if every Jordan domain $\Omega\subset X$ satisfies $$\hm^2(\Omega)\leq\frac{1}{4\pi}\cdot\length(\partial \Omega)^2.$$
\ec

In particular, spaces as in the corollary are Lipschitz $1$-connected. For spaces satisfying \eqref{eq:geom-isop-intro} with $C>\frac{1}{4\pi}$, it is in general not easy to construct non-trivial maps with suitable metric or analytic properties, as the following open problem shows:

\begin{quest}\label{quest:Lip-one-connected}
Let $X$ be a geodesic, Ahlfors $2$-regular, linearly locally connected metric space homeomorphic to $\overline{D}$. Is it true that every Lipschitz map $c\colon S^1\to X$ extends to a Lipschitz map on $\overline{D}$? If so, is $X$ Lipschitz $1$-connected?
\end{quest}

Notice that Theorem~\ref{thm:isop-ineq-for-geod-discs} in particular asserts the existence of many non-trivial Sobolev maps into a space satisfying \eqref{eq:geom-isop-intro}. Using this theorem together with the results in \cite{LWY16}, we can give a partial answer to the question above. Let $(X,d)$ be a metric space, $A\subset\R^2$, and $\alpha>0$. A map $\varphi\colon A\to X$ is said to be $(L, \alpha)$-H\"older if $d(\varphi(a), \varphi(a')) \leq L\cdot |a-a'|^\alpha$ for all $a,a'\in A$. The space $X$ is said to have the planar $\alpha$-H\"older extension property if there exists $\lambda\geq 1$ such that any $(L, \alpha)$-H\"older map $\varphi\colon A\to X$ with $A\subset \R^2$ extends to an $(\lambda L, \alpha)$-H\"older map $\bar{\varphi}\colon \R^2\to X$.  Theorem~\ref{thm:isop-ineq-for-geod-discs} together with \cite[Theorems 7.1 and 6.4]{LWY16} implies:

\bc
 Let $X$ be a complete, geodesic metric space homeomorphic to a $2$-dimensional manifold. Suppose there exists $C>0$ such that every Jordan curve in $X$ bounds a Jordan domain $U\subset X$ with 
\begin{equation*}
\hm^2(U)\leq C\cdot \length(\partial U)^2.
\end{equation*} 
Then $X$ has the planar $\alpha$-H\"older extension property for every $\alpha\in(0,1)$. 
\ec

The following special case provides our almost answer to Question~\ref{quest:Lip-one-connected}.

\bc\label{cor:Hoelder-1-connected-Ahlfors}
 Let $X$ be a geodesic, Ahlfors $2$-regular, linearly locally connected metric space homeomorphic to $\overline{D}$ or $S^2$. Then for every $\alpha\in(0,1)$ there exists $\lambda\geq 1$ such that if $\varphi\colon S^1\to X$ is an $L$-Lipschitz map then $\varphi$ extends to an $(\lambda L, \alpha)$-H\"older map on all of $\overline{D}$. 
\ec

%Recall that $X$ is said to be $\alpha$-H\"older $1$-connected if there exists $\lambda\geq 1$ such that every $L$-Lipschitz map $\varphi\colon S^1\to X$ admits an $(\lambda L, \alpha)$-H\"older extension to all of $\overline{D}$. If $\alpha=1$ then this condition is known as Lipschitz $1$-connectedness.
%
%Notice that by using the fact that quasisymmetric homeomorphisms are $\beta$-H\"older for some $\beta$ one can obtain the corollary only for sufficiently small $\alpha$.

The structure of the paper is as follows. In Section~\ref{sec:prelims} we gather definitions and some basic results which will be used throughout the paper. Section~\ref{sec:Sobolev-defs} contains background from the theory of metric space valued Sobolev maps. We furthermore summarize the main existence and regularity results for energy minimizers in metric spaces which will be used in the proof of our main theorem.
In Section~\ref{sec:top-props-Plateau} we prove Theorem~\ref{thm:Plateau-solution-cell-like-intro}. Section~\ref{sec:isop-inequality-surfaces} is devoted to the proof of Theorem~\ref{thm:isop-ineq-for-geod-discs}. We furthermore obtain Theorem~\ref{thm:almost-qc-param} which gives an almost parametrization result for surfaces as in Theorem~\ref{thm:isop-ineq-for-geod-discs}. In Section~\ref{sec:proofs-main-thms}, we prove Theorem~\ref{thm:main-thm-qs-discs-summary} and an analogous result for spheres, see Theorem~\ref{thm:BK-canonical-version}, yielding in particular the Bonk--Kleiner theorem.

\bigskip

{\bf Acknowledgements:} The authors wish to thank Katrin F\"assler, Anton Petrunin, Pekka Koskela, Kai Rajala, and Stephan Stadler for discussions. Parts of this paper were written while the second author was visiting the Scuola Normale Superiore di Pisa. He wishes to thank the institute for the hospitality and the inspiring atmosphere he enjoyed during his visit. 

\section{Preliminaries}\label{sec:prelims}

\subsection{Basic definitions and notations}
The Euclidean norm of a vector $v\in\R^2$ is denoted by $|v|$, the open unit disc in $\R^2$ by $$D:= \left\{z\in\R^2: |z|<1\right\},$$ and its closure by $\overline{D}$. 
Let $(X,d)$ be a metric space. The open ball in $X$ centered at some point $x_0$ of radius $r>0$ is denoted by $$B(x_0, r) = B_X(x_0,r) = \{x\in X: d(x,x_0)<r\}.$$
A Jordan curve in $X$ is a subset of $X$ which is homeomorphic to $S^1$. Let $\Gamma\subset X$ be a Jordan curve. A continuous map $c\colon S^1\to \Gamma$ is called a weakly monotone parametrization of $\Gamma$ if $c$ is the uniform limit of homeomorphisms $c_i\colon S^1\to \Gamma$. The length of a curve $c$ in $X$ is denoted by $\length_X(c)$ or simply by $\length(c)$.
A curve $c\colon[a,b]\to X$ is called geodesic if $\length(c) = d(c(a),c(b))$. The metric space $X$ is called geodesic if any pair of points $x,y\in X$ can be joined by a geodesic. If $c\colon(a,b)\to X$ is absolutely continuous then $c$ is metrically differentiable at almost every $t\in(a,b)$, that is, the limit $$|c'(t)|:= \lim_{s\to t}\frac{d(c(s), c(t))}{|t-s|}$$ exists, see \cite{Kir94}, \cite{KS93}.

Given $m\geq 0$ the $m$-dimensional Hausdorff measure on $X$ is denoted by $\hm^m_X$ or simply by $\hm^m$ if there is no danger of ambiguity. The normalizing constant in the definition of $\hm^m$ is chosen in such a way that $\hm^m$ coincides with the Lebesgue measure on Euclidean $\R^m$.

The following elementary observation will be used repeatedly troughout the text.

\bl\label{lem:change-metric}
 Let $(X, d)$ be a geodesic metric space homeomorphic to $D$ or $\overline{D}$ and let $\Omega\subset X$ be a Jordan domain such that $\length(\partial \Omega)<\infty$. Then the length metric $d_{\overline{\Omega}}$ on $\overline{\Omega}$ is finite for any pair of points and has the following properties:
 \begin{enumerate}
  \item The metrics $d$ and $d_{\overline{\Omega}}$ induce the same topology on $\overline{\Omega}$.
  \item For every curve $c$ in $\overline{\Omega}$ we have $\length_d(c) = \length_{d_{\overline{\Omega}}}(c)$.
  \item For every Borel set $A\subset\overline{\Omega}$ we have $\hm^2_d(A) = \hm^2_{d_{\overline{\Omega}}}(A)$.
 \end{enumerate}
 In particular, the space $(\overline{\Omega}, d_{\overline{\Omega}})$ is a geodesic metric space homeomorphic to $\overline{D}$. 
\el

The proof is a straight-forward verification and is left to the reader.

\subsection{Quasisymmetric homeomorphisms and conformal modulus}

We collect basic definitions concerning quasisymmetric mappings between metric spaces and the modulus of families of curves. We refer to \cite{TV80}, \cite{Hei01}, \cite{HK98}, \cite{HKST15} for more details.

\bd
A metric space $X$ is called linearly locally connected if there exists $\lambda\geq 1$ such that for every $x\in X$ and for all $r>0$, every pair of points in $B(x,r)$ can be joined by a continuum in $B(x,\lambda r)$ and every pair of points in $X\setminus B(x,r)$ can be joined by a continuum in $X\setminus B(x, r/\lambda)$. 
\ed

If $X$ is geodesic then the first condition in the definition of linear local connectedness is automatically satisfied. By \cite[Lemma 2.5]{BK02}, a metric space $X$ homeomorphic to a closed $2$-dimensional manifold is linearly locally connected if and only if $X$ is linearly locally contractible: there exists $\lambda\geq 1$ such that every ball $B(x,r)\subset X$ with $0<r<\lambda^{-1}\cdot\diam X$ is contractible in $B(x, \lambda r)$.

\bd
A homeomorphism $\varphi\colon M\to X$ between metric spaces $M$ and $X$ is said to be quasisymmetric if there exists a homeomorphism $\eta\colon [0,\infty)\to[0,\infty)$ such that  
\begin{equation*}
 d(\varphi(z), \varphi(a)) \leq \eta(t)\cdot d(\varphi(z), \varphi(b)).
\end{equation*}
for all $z,a,b\in M$ with $d(z,a)\leq t\cdot d(z,b)$.
\ed

Quasisymmetric homeomorphisms preserve the doubling property and linear local connectedness. Recall that a metric space $X$ is called doubling if there exists $N\geq 1$ such that every ball of radius $2r$ in $X$ can be covered by at most $N$ balls of radius $r$. 
Subsets of Ahlfors regular spaces are, in particular, doubling.

\bd
A metric space $X$ is called Ahlfors $2$-regular if there exists $L>0$ such that 
\begin{equation*}
 L^{-1}\cdot r^2\leq \hm^2(B(x,r))\leq L\cdot r^2
\end{equation*}
for all $x\in X$ and $0<r<\diam X$. 
\ed

Let $X$ be a metric space and $\Gamma$ a family of curves in $X$. A Borel function $\rho\colon X\to [0,\infty]$ is said to be admissible for $\Gamma$ if $\int_\gamma \rho\geq 1$ for every locally rectifiable curve $\gamma\in\Gamma$. We refer to \cite{HKST15} for the definition of the path integral $\int_\gamma\rho$.
The modulus of $\Gamma$ is defined by $$\MOD(\Gamma):= \inf_\rho\int_X\rho^2\,d\hm^2,$$ where the infimum is taken over all admissible functions for $\Gamma$. Note that throughout this paper, the reference measure on $X$ will always be the $2$-dimensional Hausdorff measure. By definition, $\MOD(\Gamma)=\infty$ if $\Gamma$ contains a constant curve. A property is said to hold for almost every curve in $\Gamma$ if it holds for every curve in $\Gamma_0$ for some family $\Gamma_0\subset \Gamma$ with $\MOD(\Gamma\setminus \Gamma_0)=0$. In the definition of $\MOD(\Gamma)$, the infimum can equivalently be taken over all weakly admissible functions, that is, Borel functions $\rho\colon X\to [0,\infty]$ such that $\int_\gamma \rho\geq 1$ for almost every every locally rectifiable curve $\gamma\in\Gamma$.

\bt\label{thm:weak-modulus-qc-implies-qs}
 Let $X$ be a metric space which is homeomorphic to $\overline{D}$ and satisfies, for some $L>0$, $$\hm^2(B(x, r))\leq L\cdot r^2$$ for every $x\in X$ and $r>0$. Let $Q\geq 1$ and suppose $u\colon \overline{D} \to X$ is a homeomorphism satisfying $$ \MOD(\Gamma)\leq Q\cdot \MOD(u\circ\Gamma)$$ for every family $\Gamma$ of curves  in $\overline{D}$. Then $u$ is quasisymmetric if and only if $X$ is doubling and linearly locally connected.
\et

Here, $u\circ\Gamma$ denotes the family of curves $u\circ \gamma$ with $\gamma\in\Gamma$. The proof of the theorem is a simple variation of the proof of \cite[Theorem 4.7]{HK98}, which we provide for completeness in the appendix. The theorem furthermore holds with $\overline{D}$ replaced by $S^2$.

\subsection{Topological preliminaries}\label{sec:prelim-top}

We recall some topological notions and results which we will need in Section~\ref{sec:top-props-Plateau}. For details we refer to \cite{HNV04}, \cite{Edw78}, \cite{Dav07}.

\bd\label{def:monotone-light}
 Let  $X$ and $Y$ be metric spaces and $v\colon X\to Y$ a continuous map. If $v^{-1}(y)$ is connected for every $y\in Y$ then $v$ is called monotone. If $v^{-1}(y)$ is totally disconnected for every $y\in Y$ then $v$ is called light.
\ed

The monotone-light factorization theorem due to Eilenberg and Whyburn asserts that for every continuous, surjective map $v \colon X\to Y$ between compact metric spaces $X$ and $Y$ there exist a compact metric space $Z$ and continuous surjective maps $v_1\colon X\to Z$ and $v_2\colon Z\to Y$ such that $v_1$ is monotone, $v_2$ is light, and $v= v_2\circ v_1$. The fibers $v_1^{-1}(z)$ are exactly the connected components of $v^{-1}(v_2(z))$. 

We will furthermore need the notion of cell-like spaces and cell-like maps.

\bd
 A compact metric space is called cell-like if it admits an embedding into the Hilbert cube $Q$ in which it is null-homotopic in every neighborhood of itself. A continuous surjection $f\colon X\to Y$ between metric spaces $X$ and $Y$  is called cell-like if $f^{-1}(q)$ is cell-like (in particular compact) for every $q\in Y$.
\ed

A compact subset $X$ of $S^2$ is cell-like if and only if $X$ is connected and $S^2\setminus X$ is non-empty and connected. A closed connected subset $X$ of $\overline{D}$ is cell-like if and only if each connected component of $\overline{D}\setminus X$ intersects $\partial D$. 
A cell-like set contained in a Jordan curve is a point or a topological arc. A compact $1$-dimensional cell-like metric space $X$ is unicoherent, see \cite[p.~332]{HNV04} and  \cite[p.~97]{HNV04}. Recall that a connected metric space $X$ is unicoherent if for all closed connected subsets $A,B\subset X$ with $X=A\cup B$ the intersection $A\cap B$ is connected.

Let $X$ and $Y$ be absolute neighborhood retracts. A continuous surjective map $f\colon X\to Y$ is cell-like if and only if for every open set $U\subset Y$ the restriction $$f|_{f^{-1}(U)}\colon f^{-1}(U)\to U$$ is a homotopy equivalence. In particular, if $f$ is cell-like then for every open connected (respectively contractible) set $U\subset Y$ the preimage $f^{-1}(U)$ is connected (respectively contractible).

\bl\label{lem:prelim-top-subsets}
 Let $W\subset S^2$ be an open, connected set and let $K\subset W$ be a compact set all of whose connected components are cell-like. Then $W\setminus K$ is connected. Moreover, for every connected component $T$ of $K$ there exist arbitrarily small neighborhoods $V\subset W$ of $T$ which are homeomorphic to $D$ and such that $\overline{V}$ is homeomorphic to $\overline{D}$ and $\partial V$ does not intersect $K$.
\el

\begin{proof}
 Let $Y$ be the space obtained from $W$ by identifying the points in each connected component of $K$, endowed with the finest topology such that the natural projection $\pi \colon W\to Y$ is continuous. Then $\pi$ is a cell-like map and hence $Y$ is homeomorphic to $W$ by Moore's theorem, see e.g.~\cite{Dav07}. Moreover, $\pi(K)$ is totally disconnected in $Y$ and hence $Y\setminus \pi(K)$ is connected. It follows that $W\setminus K$ is also connected. This proves the first statement of the lemma.
 
Let now $T\subset K$ be a connected component of $K$. Then there exist arbitrarily small neighborhoods $U$ of the point $\pi(T) \in Y$ such that $U$ is homeomorphic to $D$,  $\overline{U}$ is homeomorphic to $\overline{D}$ and such that the circle $\partial U$ does not intersect the totally disconnected set $\pi(K)$.  Then $V=\pi^{-1} (U)$ is an arbitrarily small neighborhood of $T$ which is homeomorphic to $D$, such that $\overline{V}$ is homeomorphic to $\overline{D}$ and such that the circle $\partial V$ does not intersect  $K$.
\end{proof}

\bp\label{prop:equiv-cell-like-monotone}
For a continuous surjective map $f\colon \overline{D}\to\overline{D}$ the following statements are equivalent:
\begin{enumerate}
 \item $f$ is monotone.
 \item $f$ is cell-like.
 \item $f$ is a uniform limit of homeomorphisms $f_i\colon\overline{D}\to\overline{D}$.
\end{enumerate}
\ep

\begin{proof}
 Clearly, (ii) implies (i). Property (i) is equivalent to (iii) by \cite{You48}. Finally, (iii) implies (ii) since the uniform limit of cell-like maps between compact absolute neighborhood retracts is cell-like by \cite[Theorem 17.4]{Dav07}.
\end{proof}

\section{Metric space valued Sobolev maps}\label{sec:Sobolev-defs}

We recall some definitions from the theory of metric space valued Sobolev mappings based on upper gradients \cite{Shan00}, \cite{HKST01}, \cite{HKST15} as well as the results concerning existence and regularity of energy minimizing discs established in \cite{LW15-Plateau}, \cite{LW16-harmonic}, \cite{LW16-energy-area}. 

Let $(X,d)$ be a complete metric space and $\Omega\subset \R^2$ a bounded domain.  A Borel function $g\colon \Omega\to [0,\infty]$ is said to be an upper gradient of a map $u\colon \Omega\to X$ if
\begin{equation}\label{eq:upper-grad}
 d(u(\gamma(a)), u(\gamma(b)))\leq \int_\gamma g
\end{equation}
 for every rectifiable curve $\gamma\colon[a,b]\to \Omega$. If \eqref{eq:upper-grad} only holds for almost every curve $\gamma$ then $g$ is called a weak upper gradient of $u$. A weak upper gradient $g$ of $u$ is called minimal weak upper gradient of $u$ if $g\in L^2(\Omega)$ and if for every weak upper gradient $h$ of $u$ in $L^2(\Omega)$ we have $g\leq h$ almost everywhere on $\Omega$.

Denote by $L^2(\Omega, X)$ the collection of measurable and essentially separably valued maps $u\colon \Omega\to X$ such that the function $u_x(z):= d(u(z), x)$ belongs to $L^2(\Omega)$ for some and thus any $x\in X$. A map $u\in L^2(\Omega, X)$ belongs to the (Newton-) Sobolev space $N^{1,2}(\Omega, X)$ if $u$ has a weak upper gradient in $L^2(\Omega)$. Every such map $u$ has a minimal weak upper gradient $g_u$, unique up to sets of measure zero, see \cite[Theorem 6.3.20]{HKST15}.
The {\it Reshetnyak energy} of a map $u\in N^{1,2}(\Omega, X)$ is defined by $$E_+^2(u):= \|g_u\|_{L^2(\Omega)}^2.$$

If $u\in N^{1,2}(\Omega, X)$ then for almost every $z\in \Omega$ there exists a unique semi-norm on $\R^2$, denoted by $\apmd u_z$ and called the approximate metric derivative of $u$, such that 
$$\ap \lim _{y\to z}  \frac {d(u(y),u(z))- \apmd u_z(y-z)} {|y-z|} =0,$$ see \cite{Kar07} and \cite[Proposition 4.3]{LW15-Plateau}. Here, $\ap\lim$ denotes the approximate limit, see \cite{EG92}. 
 
The following notion of parametrized area was introduced in \cite{LW15-Plateau}.

\bd
 The parametrized (Hausdorff) area of a map $u\in N^{1,2}(\Omega, X)$ is defined by $$\Area(u)= \int_\Omega \jac(\apmd u_z)\,dz,$$ where the Jacobian $\jac(s)$ of a semi-norm $s$ on $\R^2$ is the Hausdorff $2$-measure on $(\R^2, s)$ of the unit square if $s$ is a norm and zero otherwise.
 \ed
 
The area of the restriction of $u$ to a measurable set $B\subset \Omega$ is defined analogously. If $u$ is injective and satisfies Lusin's condition (N) then $\Area(u) = \hm^2(u(\Omega))$ by the area formula \cite{Kir94}, \cite{Kar07}.

Recall that, by John's theorem, the unit ball with respect to a norm $\|\cdot \|$ on $\R^2$ contains a unique ellipse of maximal area, called John's ellipse of $\|\cdot\|$. We will need the following definition from \cite{LW16-energy-area}. 

\bd\label{def:inf-isotropic}
 A map $u\in N^{1,2}(\Omega, X)$ is called infinitesimally isotropic if for almost every $z\in \Omega$ the semi-norm $\apmd u_z$ is either zero or is a norm and the John ellipse of $\apmd u_z$ is a Euclidean disc.
\ed

We call a map $u\in N^{1,2}(\Omega, X)$ {\it infinitesimally $Q$-quasiconformal} if
\begin{equation}\label{eq:inf-qc}
(g_u(z))^2 \leq Q\cdot \jac(\apmd u_z)
\end{equation}
for almost every $z\in \Omega$. If $u\in N^{1,2}(\Omega, X)$ is infinitesimally isotropic then it is infinitesimally $Q$-quasi\-conformal with $Q=\frac{4}{\pi}$, see \cite{LW16-energy-area}. 

%By \cite[Lemma 7.2]{LW15-Plateau}, for every $u\in N^{1,2}(D, X)$ we have $$\jac(\apmd u_z)\leq (g_u(z))^2$$ for almost every $z\in D$. 

If $u\in N^{1,2}(D, X)$ then for almost every $v\in S^1$ the curve $t\mapsto u(tv)$ with $t\in[1/2, 1)$ is absolutely continuous. The trace of $u$ is defined by $$\trace(u)(v):= \lim_{t\nearrow 1}u(tv)$$ for almost every $v\in S^1$. It can be shown that $\trace(u)\in L^2(S^1, X)$, see \cite{KS93}. If $u$ is the restriction to $D$ of a continuous map $\hat{u}$ on $\overline{D}$ then $\trace(u)=\hat{u}|_{S^1}$.

Given a Jordan curve $\Gamma\subset X$ we denote by $\Lambda(\Gamma, X)$ the possibly empty family of maps $v\in N^{1,2}(D, X)$ whose trace has a continuous representative which weakly monotonically parametrizes $\Gamma$. The following theorem summarizes the existence and regularity properties of energy minimizers which we will need in this article and which were proved in \cite{LW15-Plateau}, \cite{LW16-energy-area}, \cite{LW16-harmonic}. Note that the results in these papers are stated using a different but equivalent definition of Sobolev mapping and Reshetnyak energy, see \cite{Res04} and \cite[Theorem 7.1.20]{HKST15}.

\bt\label{thm:existence-regularity-energy-min}
 Let $X$ be a proper metric space admitting a quadratic isoperimetric inequality. Let $\Gamma\subset X$ be a Jordan curve of finite length. Then $\Lambda(\Gamma, X)$ is non-empty and contains an element $u$ satisfying $$E_+^2(u) = \inf\left\{E_+^2(v): v\in\Lambda(\Gamma, X)\right\}.$$ Any such $u$ is infinitesimally isotropic and has a representative which is continuous on $D$ and extends continuously to $\overline{D}$.
\et

\begin{proof}
The existence of an energy minimizer in $\Lambda(\Gamma, X)$ follows from \cite[Theorem 7.6]{LW15-Plateau}. Continuity of energy minimizers up to the boundary is a consequence of \cite[Theorem 4.4]{LW16-energy-area} or \cite[Theorem 1.3]{LW16-harmonic}. Infinitesimal isotropy follows from \cite[Lemmas 3.2 and 4.1]{LW16-energy-area}, see also \cite[Lemma 6.5]{LW15-Plateau}. 
\end{proof}

We end this section with the following proposition. See \cite[Theorem 1.1]{Wil12} for an analogous result.

\bp\label{prop:ae-qc-implies-weak-modulus-qc}
 Let $X$ be a complete metric space and $u\colon \overline{D}\to X$ continuous and monotone. If $u\in N^{1,2}(D, X)$ and $u$ is infinitesimally $Q$-quasiconformal 
 then 
 \begin{equation*}
 \MOD(\Gamma)\leq Q\cdot \MOD(u\circ\Gamma)
 \end{equation*}
for every family $\Gamma$ of curves in $\overline{D}$.
\ep

\begin{proof}
Let $g_u$ be the mininal weak upper gradient of $u$ on $D$. We first claim that the upper gradient inequality \eqref{eq:upper-grad} holds with $g=g_u$ for almost every rectifiable curve $\gamma$ in $\overline{D}$ instead of $D$. Indeed, $u$ extends to a Newton-Sobolev map on the open disc $B(0,2)\subset\R^2$ and hence $u\circ\gamma$ is absolutely continuous for almost every rectifiable curve in $\overline{D}$, parametrized by arc-length. Since almost every curve $\gamma$ in $\overline{D}$ intersects the boundary $S^1$ in a set of Hausdorff $1$-measure zero, the claim now follows from \eqref{eq:upper-grad}. 

Let $\Gamma$ be a family of curves in $\overline{D}$. Then for almost every rectifiable curve $\gamma\in\Gamma$, parametrized by arc-length, the curve $u\circ\gamma$ is absolutely continuous, and for almost every $t$ we have $\gamma(t)\not\in S^1$ and $u\circ\gamma$ is metrically differentiable at $t$ with
\begin{equation}\label{eq:approx-metr-deriv-composition}
 |(u\circ\gamma)'(t)| \leq g_u(\gamma(t)),
\end{equation}
see \cite[Proposition 6.3.3]{HKST15}.
Let $\varrho\colon X\to [0,\infty]$ be an admissible function for the family $u\circ\Gamma$ and define $\bar{\varrho}: = g_u\cdot (\varrho\circ u)$ on $D$. Then for every curve $\gamma\in\Gamma$ with the properties above, inequality \eqref{eq:approx-metr-deriv-composition} yields $$\int_\gamma \bar{\varrho} = \int g_u(\gamma(t)) \cdot \varrho(u\circ\gamma(t))\,dt \geq \int |(u\circ\gamma)'(t)|\cdot\varrho(u\circ\gamma(t))\,dt = \int_{u\circ\gamma}\varrho\geq 1,$$ so $\bar{\varrho}$ is weakly admissible for $\Gamma$.

By \cite[Proposition 3.2]{LW15-Plateau}, there exists a set $A\subset D$ of measure zero such that the restriction $u|_{D\setminus A}$ has Lusin's property (N). For $x\in u(D)$ let $N(u,x)$ denote the number of points $z\in D$ with $u(z) = x$. Clearly, $N(u,x)$ equals $1$ or $\infty$ for every $x\in u(D)$ because $u$ is monotone. Since $$\int_{u(D\setminus A)} N(u, x)\,d\hm^2(x) = \int_{D\setminus A}\jac(\apmd u_z)\,dz = \Area(u)<\infty$$ by the area formula \cite{Kir94}, \cite{Kar07} it follows that $N(u,x)= 1$ for $\hm^2$-almost every $x\in u(D\setminus A)$.
Thus, since $u$ is infinitesimally $Q$-quasi\-con\-for\-mal and monotone it follows again from the area formula that
 \begin{equation*}
   \int_D \bar{\varrho}(z)^2\,dz \leq Q\cdot \int_D \varrho\circ u(z)^2\cdot \jac(\apmd u_z)\,dz \leq Q\cdot \int_X\varrho(x)^2\,d\hm^2(x).
  \end{equation*}
Since $\varrho$ was an arbitrary admissible function for $u\circ\Gamma$ it therefore follows that $$\MOD(\Gamma)\leq Q\cdot \MOD(u\circ\Gamma).$$ This completes the proof.
\end{proof}

\section{Topological properties of energy minimizers}\label{sec:top-props-Plateau}

In this section we prove Theorem~\ref{thm:Plateau-solution-cell-like-intro}. The key ingredient is the following topological result.

\bt\label{thm:cell-like-homeo}
Let $X$ be a geodesic metric space homeomorphic to $\overline{D}$. Let $v\colon\overline{D}\to X$ be a continuous and surjective map satisfying the following properties:
\begin{enumerate}
 \item The restriction of $v$ to $S^1$ is a weakly monotone parametrization of $\partial X$. 
 \item Whenever $T \subset X$ is a single point or biLipschitz homeomorphic to a closed interval then every connected component of $v^{-1} (T)$ is cell-like.
\end{enumerate}
Then $v$ is a cell-like map.
\et

Before proving Theorem~\ref{thm:cell-like-homeo} we first provide:

\begin{proof}[Proof of Theorem~\ref{thm:Plateau-solution-cell-like-intro}]
Let $u$ be as in the theorem. Then $u$ is infinitesimally quasiconformal by Theorem~\ref{thm:existence-regularity-energy-min}. Moreover, $u$ minimizes the inscribed Riemannian area $\Area_{\mu^i}$ among all maps in $\Lambda(\partial X, X)$ by \cite[Theorem 4.3 and Corollary 3.3]{LW16-energy-area}.

Due to Proposition~\ref{prop:equiv-cell-like-monotone}, in order to prove the theorem it is enough to show that $u$ satisfies the hypotheses of Theorem~\ref{thm:cell-like-homeo}. Clearly, the map $u$ is surjective and satisfies property (i). In order to see that $u$ also satisfies property (ii) we argue by contradiction and assume that there exists $T$ as in (ii) such that some connected component $K$ of $u^{-1}(T)$ is not cell-like. Thus, there exists a connected component of $\overline{D} \setminus K$ which does not intersect $S^1$. In particular, there also exists a connected component $U$ of $\overline{D}\setminus u^{-1} (T)$ which does not intersect $S^1$.
Since $T$ is an absolute Lipschitz retract there exists a Lipschitz projection $P\colon X\to T$. Define a map by $w:= P\circ u$.  Then the restrictions $w|_U$ and $u|_U$ have the
same trace in the sense of \cite[Definition 4.1]{LW16-intrinsic}. Hence, by  \cite[Lemma 4.2]{LW16-intrinsic}, we may replace $u|_U$ by $w|_U$ and obtain another map $u_1\in N^{1,2}(D, X)$ with the same trace as $u$ and, in particular, $u_1\in \Lambda(\partial X, X)$. 
Since $\hm^2(T)=0$ the inscribed Riemannian area of $w$ is zero. Since $u$ minimizes the inscribed Riemannian area  it follows that the inscribed Riemannian area of $u|_U$ is zero. Since $u$ is infinitesimally quasiconformal the Reshetnyak energy of $u|_U$ must be zero as well.
Therefore $u|_U$ is constant and hence $u(U)$ is contained in $T$, a contradiction. Thus every connected component of $u^{-1}(T)$ is cell-like and hence $u$ satisfies (ii). This completes the proof.
\end{proof}

The rest of this section is devoted to the proof of Theorem~\ref{thm:cell-like-homeo}. Let $v$ be as in the statement of the theorem. We first claim that it is enough to consider the case that $v$ is, in addition, a light map. Indeed, by the monotone-light factorization theorem there exist a metric space $Z$ and continuous, surjective maps $v_1\colon \overline{D}\to Z$ and $v_2\colon Z\to X$ such that $v_1$ is monotone, $v_2$ is light, and $v$ factors as $v = v_2\circ v_1$. Moreover, the fibers $v_1^{-1}(z)$ are exactly the connected components of $v^{-1}(v_2(z))$. It thus follows from the properties of $v$ that $v_1$ and the restriction $v_1|_{S^1}$ are cell-like maps. Consequently, $Z$ is homeomorphic to $\overline{D}$ and $v_1$ is the uniform limit of homeomorphisms, see \cite[Corollary 7.12]{LW16-intrinsic}. We now identify $Z$ with $\overline{D}$ via a homeomorphism and show:

\bl\label{lem:v2-has-props}
 The map $v_2$ satisfies properties (i) and (ii) of Theorem~\ref{thm:cell-like-homeo}.
\el

\begin{proof}
 We first prove that $v_2$ has property (i). As mentioned above, we identify $Z$ with $\overline{D}$ via a homeomorphism. Since $v_1$ is the uniform limit of homeomorphisms it follows that $v_1(S^1) = S^1$ and hence $v_2(S^1) = \partial X$. Let $x\in\partial X$. We must show that $v_2^{-1}(x)\cap S^1$ consists of a single point. Let $z, z'\in S^1$ be such that $v_2(z)=x=v_2(z')$. The preimages $v_1^{-1}(z)$ and $v_1^{-1}(z')$ are connected components of $v^{-1}(x)$, both having non-trivial intersection with $S^1$. Since $v^{-1}(x)\cap S^1$ is an interval or a point we must therefore have $z=z'$. This proves property (i).
 
As for property (ii) let $T\subset X$ be a point or biLipschitz homeomorphic to a closed interval and let $K\subset \overline{D}$ be a connected component of $v_2^{-1}(T)$. We must show that $K$ is cell-like. Since $v_1$ is surjective and monotone we have that $K':= v_1^{-1}(K)$ is connected and thus $K'$ is a connected component of $v^{-1}(T)$. In particular, $K'$ is cell-like and hence $K$ is cell-like as well by \cite[Theorem 1.4]{Lac69}. This proves that $v_2$ satisfies property (ii) and completes the proof of the lemma.
\end{proof}

Now, if $v_2$ is cell-like then so is $v$ because $v_1$ is cell-like. This together with Lemma~\ref{lem:v2-has-props} shows our claim. Thus, for the proof of Theorem~\ref{thm:cell-like-homeo} it is indeed enough to consider only the case that $v$ is, in addition, a light map. We henceforth assume that $v$ is as in the statement of Theorem~\ref{thm:cell-like-homeo} and that $v$ is also light. We must show that $v$ is injective. For this, we first prove some auxiliary results.

\bl\label{lem:bilipap}
Let $X$ be a geodesic metric space and let $\Gamma\subset  X$ be a topological arc connecting two points $a,b\in X$.
Then for every $\varepsilon>0$ there exists a biLipschitz curve contained in the $\varepsilon$-neighborhood of $\Gamma$ and connecting $a$ and $b$.
\el

A similar statement holds for Jordan curves.

\begin{proof}
Let $\varepsilon>0$ and choose a piecewise geodesic curve $\gamma$ which is  contained in the $\frac{\varepsilon}{2}$-neighborhood of $\Gamma$ and connects the endpoints $a$ and $b$ of $\Gamma$.
By choosing an appropriate subcurve of $\gamma$ we may assume that $\gamma$ is an injective piecewise biLip\-schitz curve. By changing the curve step by step near its vertices it suffices to prove the following claim.

\emph{Claim:}  let $s>0$ and let $\eta\colon [-s,s] \to X$ be an injective curve such that the restrictions $\eta|_{[0,s]}$ and $\eta|_{[-s,0]}$ are geodesics  parametrized by arc-length.
Then there exist arbitarily small $t\in(0,s)$ such that after replacing $\eta|_{[-t,t]}$ by a geodesic from $\eta(-t)$ to $\eta(t)$ we obtain a biLipschitz curve.

To prove the claim, note first that the Lipschitz function $f(t)=d(\eta(-t),\eta (t))$ satisfies $f(0)=0$ and is strictly positive for $t>0$. Thus we find arbitrarily small
$t >0$  for which $f'(t)$ exists and is strictly positive.  Fix such $t$ and set $\delta:= f'(t)$. Choose a geodesic $c_t$ from $\eta(-t)$ to $\eta(t)$, parametrized by arc-length and such that $c_t(0)= \eta(-t)$. The triangle inequality and the fact that $f'(t)=\delta>0$ yield that, for all sufficiently small $r>0$, we have $$d(\eta(-t-r), c_t(r))\geq \frac{\delta \cdot r}{2}.$$ Since this holds for every such geodesic $c_t$ it follows that $c_t$ can only intersect $\eta|_{[-s,-t]}$ at $\eta(-t)$. Moreover, the inequality above together with the triangle inequality imply that the concatenation of $\eta|_{[-s,-t]}$ with $c_t$ is a biLipschitz curve locally around $\eta(-t)$ and hence also globally. The same argument shows that the concatenation of $c_t$ with $\eta|_{[t,s]}$ is a biLipschitz curve, hence the concatenation of all three curves is biLipschitz. This proves the claim and completes the proof of the lemma.
\end{proof}

\bl\label{lem:unique-preimage-on-S1}
For every $x\in \partial X$ the preimage $v^{-1}(x)$ consists of exactly one point.
\el

\begin{proof}
Since $v$ is light and its restriction to $S^1$ is a weakly monotone parametrization of $\partial X$, the point $x$ has exactly one preimage in $S^1$. We must show that $x$ has no preimages in $D$. We argue by contradiction and assume that there exists $z\in D$ such that $v(z)=x$. There exists an open neighborhood $U\subset D$ of  $z$ homeomorphic to $D$ such that $\overline{U}$ is homeomorphic to $\overline{D}$ and $\partial U$ does not intersect the totally disconnected set $v^{-1}(x)$.  In an arbitrary small neighborhood of $x$ we find a simple arc $S$ connecting two points on $\partial X$ on different sides of $x$ such that $S$ does not intersect $v(\partial U)$.  If $S$ is sufficiently close to $x$ then $S$ separates $v(U)$. By Lemma~\ref{lem:bilipap} we may assume that $S$ is a biLipschitz curve. Since $S$ separates $v(U)$ the set $v^{-1}(S)$ must separate $U$. However, $v^{-1} (S)\cap \partial  U =\emptyset$ and any connected component of $v^{-1} (S)$ is cell-like by assumption. Hence the set $K:= v^{-1}(S)\cap U$ is compact and all its connected components are cell-like. So, $U\setminus v^{-1} (S)$ is connected by Lemma~\ref{lem:prelim-top-subsets}. However, this contradicts the fact that $v^{-1}(S)$ separates $U$ and completes the proof.
\end{proof}

\bl\label{lem:image-of-preimage}
Let $T\subset X$ be biLipschitz homeomorphic to a closed interval. Then $v(C)=T$ for every connected component $C$ of $v^{-1} (T)$.
\el

\begin{proof}
 We argue by contradiction and assume that there exists a conntected component $C$ of $v^{-1}(T)$ such that $v(C)\not= T$ and thus $v(C)$ is a compact subarc of $T$.  After possibly replacing $T$ by a non-trivial subarc which intersects $v(C)$ exactly in one
point, and replacing $C$  by a connected component in $C$ of the preimage of this subarc we may assume that $v(C)$ is just one endpoint of $T$, which we call $p$. Since $v$ is a light map it follows that $C$ consists of a single point $z\in \overline{D}$.

 We first assume that $z\not\in S^1$. By Lemma~\ref{lem:prelim-top-subsets} there exists an arbitrarily small open neighborhood $U\subset D$ of $z$ homeomorphic to $D$ such that $\overline{U}$ is homeomorphic to $\overline{D}$ and $\partial U$ does not intersect $v^{-1} (T)$. Let $T'$ be a non-trivial compact subarc of $T\setminus \{p\}$. Choosing $U$ sufficiently small we may assume that $T'$ has no preimage in $\overline{U}$. 
 
 We find an arc $S$ in an arbitrary small neighborhood of $T$ which connects two different points on $T'$ and which together with the corresponding part of $T'$ defines a Jordan curve  $\Gamma$ which separates $v(U)$.
 Choosing $S$ sufficiently close to $T$ we may assume that $S$ does not intersect $v (\partial U)$. Using Lemma~\ref{lem:bilipap} we may furthermore assume that $S$ is a biLipschitz curve. Since $T'$ does not have any preimage in $\overline{U}$ it follows that the preimage of $\Gamma$ in $\overline{U}$  coincides with the preimage of $S$ in $\overline{U}$. Moreover, this preimage does not intersect $\partial U$.  Thus, the set $K:= v^{-1} (\Gamma)\cap U$ is compact and every connected component of $K$ is cell-like and hence $U\setminus v^{-1}(\Gamma)$ is connected by Lemma~\ref{lem:prelim-top-subsets}. This however is impossible since $\Gamma$ separates $v(U)$. This contradiction finishes the proof in the case that $z\not\in S^1$. 
 
The proof in the case $z\in S^1$ is analogous and is left to the reader.
\end{proof}

Using the lemmas above we now show that $v$ is injective, which will complete the proof of Theorem~\ref{thm:cell-like-homeo}. Let $x\in X$. If $x\in \partial X$ then $v^{-1}(x)$ consists of a single point by Lemma~\ref{lem:unique-preimage-on-S1}, so we may assume that $x\not\in\partial X$. We first suppose that $x$ is contained in a biLipschitz arc $T\subset X$ with endpoints $x_{\pm}\in\partial X$. By Lemma~\ref{lem:unique-preimage-on-S1}, the two points $x_{\pm}$ have unique preimage points $z_{\pm}\in S^1$. By Lemma~\ref{lem:image-of-preimage}, any connected component of $v^{-1} (T)$ must contain both points $z_{\pm}$ and hence the set $\Gamma=v^{-1} (T)$ is connected and thus cell-like. 
Since $v$ is a light map, the set $\Gamma$ is $1$-dimensional, see \cite[p.~311]{HNV04}, and hence $\Gamma$ is unicoherent. Now, let $T^+$ and $T^-$ be the subintervals into which $x$ subdivides $T$. As above the preimages $\Gamma^{\pm}$  of $T^{\pm}$ must be connected. Since $\Gamma$ is unicoherent the intersection $\Gamma ^+ \cap \Gamma ^-$ is connected. Since this intersection is exactly the totally disconnected fiber $v^{-1}(x)$ it follows that $v^{-1}(x)$ has exactly one point. 

Now, let $x\in X\setminus \partial X$ be arbitrary. Connect $x$ by a geodesic $S$ with a point on $\partial X$. The construction used in Lemma~\ref{lem:bilipap} shows that any part of $S$ which does not contain $x$ can be extended to a biLipschitz arc $T$ connecting two points on $\partial X$. Thus for any point $y$  on $S\setminus \{x \}$ the preimage of $y$ contains only one point.   Since $v^{-1}(x)$ is totally disconnected we deduce from Lemma~\ref{lem:image-of-preimage} that  $v^{-1} (x)$ has also only one point. This shows that $v$ is injective and completes the proof of Theorem~\ref{thm:cell-like-homeo}.

\section{A quadratic isoperimetric inequality for metric surfaces}\label{sec:isop-inequality-surfaces}

In this section we prove Theorem~\ref{thm:isop-ineq-for-geod-discs} from the introduction as well as Corollary~\ref{cor:linloccon-isop-Sob} below. In Theorem~\ref{thm:almost-qc-param}, we moreover obtain an almost parametrization result for spaces as in Theorem~\ref{thm:isop-ineq-for-geod-discs}. We begin with the following result which proves the first part of Theorem~\ref{thm:isop-ineq-for-geod-discs}.

\bt\label{thm:isop-ineq-for-geod-discs-not-optimal}
 Let $X$ be a complete, geodesic metric space homeomorphic to a $2$-dimensional manifold. Suppose there exists $C$ such that every Jordan curve $\Gamma\subset X$ bounds a Jordan domain $\Omega\subset X$ with 
\begin{equation*}
 \hm^2(\Omega)\leq C\cdot \length(\Gamma)^2.
\end{equation*}
Then $X$ admits a quadratic isoperimetric inequality with isoperimetric constant only depending on $C$.  
\et

We need some preparations. Recall that a metric space $Y$ is said to be $L$-Lipschitz $1$-connected up to some scale if there exists $\lambda_0>0$ such that every $\lambda$-Lipschitz curve $c\colon S^1\to Y$ with $\lambda\leq \lambda_0$ extends to an $L\lambda$-Lipschitz map $\varphi\colon \overline{D}\to Y$.

\bp\label{prop:fill-Jordan-curve}
 Let $X$ be as in Theorem~\ref{thm:isop-ineq-for-geod-discs-not-optimal} and suppose $Y$ is a metric space which contains $X$ and which is $L$-Lipschitz $1$-connected up to some scale. Then  every injective Lipschitz curve $c\colon S^1\to X$ extends to a Lipschitz map $\varphi\colon\overline{D}\to Y$ with $$\Area(\varphi)\leq C'\cdot \length(c)^2$$ for some constant $C'$ only depending on $C$ and $L$.
\ep

Notice that we do not impose a bound on the length of the curve $c$ and that the constant $C'$ is independent of the scale up to which $Y$ is Lipschitz $1$-connected.

\begin{proof}
 Let $c\colon S^1\to X$ be an injective Lipschitz curve. Since $c$ is homotopic to its constant speed parametrization via a Lipschitz homotopy of zero area, see \cite[Lemma 3.6]{LWY16}, we may assume that $c$ is parametrized proportional to arc-length. Let $\Omega\subset X$ be the Jordan domain of smallest area and boundary $c$. It follows that $\hm^2(U)\leq C\cdot\length(\partial U)^2$ for every Jordan domain $U\subset \Omega$.
 
Denote by $d_{\overline{\Omega}}$ the length metric on $\overline{\Omega}$ and set $Z:= (\overline{\Omega}, d_{\overline{\Omega}})$. By Lemma~\ref{lem:change-metric}, the space $Z$ is geodesic and homeomorphic to $\overline{D}$. Moreover, the length of the boundary circle $\partial Z$ as well as the Hausdorff $2$-measure of $Z$ are finite. Finally, for every Jordan domain $U\subset Z$ we have $\hm_Z^2(U)\leq C\cdot \length_Z(\partial U)^2$.

Let $n\in\N$ be sufficiently large, to be determined later. By \cite[Theorem 4.1]{LWY16}, there exists a triangulation $\tau$ of $Z$ consisting of at most $K\cdot n^2$ triangles of $d_{\overline{\Omega}}$-diameter at most $\frac{\length(c)}{n}$ each, and such that every edge contained in $\partial Z$ has length at most $\frac{\length(c)}{n}$. Here, $K$ only depends on $C$. By a triangulation of $Z$ we mean a homeomorphism from $Z$ to a combinatorial $2$-complex $\tau$ in which every $2$-cell is a triangle. We endow $Z$ with the induced cell structure from $\tau$. For $i=0,1,2$ the $i$-skeleton of $\tau$ will be denoted $\tau^{(i)}$ and may thus be viewed as a subset of $Z$. The $2$-cells of $\tau$ will also be called triangles in $\tau$. 
 
Put on each triangle in $\tau$ a piecewise Euclidean metric which makes it an equilateral Euclidean triangle of sidelength one. Let $d_\tau$ denote the resulting length metric on $\tau$. Then $\Sigma:= (\tau, d_\tau)$ is biLipschitz homeomorphic to $\overline{D}$.  Denote by $P\colon Z\to X$ the identity map and observe that $P$ is $1$-Lipschitz and a homeomorphism onto $\overline{\Omega}$. 

We construct a Lipschitz map $\psi\colon \Sigma \to Y$ as follows. For each $v\in \tau^{(0)}$ we set $\psi(v):= P(v)$, where we have identified $v$ with its image in $Z$. Let now $e=[v,w]\in \tau^{(1)}$ be an edge. If $e$ is contained in the boundary circle $\partial \Sigma$ then we let $\psi|_{e}$ be a constant speed parametrization of the part of $c$ between the points $\psi(v)$ and $\psi(w)$. Otherwise, let $\psi|_{e}$ be a constant speed geodesic from $\psi(v)$ to $\psi(w)$ in $X$. It follows that, for every triangle $F\in \tau^{(2)}$, the map $\psi|_{\partial F}$ is $M \length(c)n^{-1}$-Lipschitz, where $M$ is a universal constant. Hence, if $n$ was chosen sufficiently large then there exists an $M'\length(c)n^{-1}$-Lipschitz extension $\psi|_F\colon F\to Y$ for some $M'$ only depending on $L$. Now, let $\varrho\colon S^1\to \partial \Sigma$ be the piecewise constant speed map such that $\psi\circ\varrho = c$. Then, $\varrho$ is biLipschitz and thus extends to a biLipschitz homeomorphism $\varrho\colon\overline{D}\to\Sigma$ by \cite[Theorem A]{Tuk80}. Thus, the map $\varphi:= \psi\circ \varrho$ is Lipschitz and satisfies $\varphi|_{S^1} = c$ and $$\Area(\varphi) = \sum_{F\in\tau^{(2)}} \Area(\psi|_F) \leq \frac{\sqrt{3}K}{4}\cdot n^2\cdot \left(M'\cdot\frac{\length(c)}{n}\right)^2 = C'\cdot\length(c)^2$$ for some constant $C'$ only depending on $C$ and $L$. This concludes the proof.
\end{proof}

The next lemma allows one to pass from non-injective Lipschitz curves to injective ones. For $m\geq 0$, let $A_m$ be the closed unit disc $\overline{D}$ with $m$ pairwise separated open Euclidean discs removed. More precisely, $$A_m= \overline{D}\setminus \cup_{j=1}^m D_j$$ for some open Euclidean discs $D_j$
such that $\overline{D}_j\subset D$ are pairwise disjoint.

\bl\label{lem:jordanize-curves}
 Let $X$ be a geodesic metric space and $c\colon S^1\to X$ a Lipschitz curve. Given points $s_0, \dots,s_N\in S^1$ there exist $m\geq 0$ and a Lipschitz map $\varphi\colon A_m\to X$ of zero area satisfying the following properties:
 \begin{enumerate}
  \item $\varphi|_{S^1}$ is a piecewise geodesic with $\varphi(s_j) = c(s_j)$ for $j=0,\dots, N$.
  \item Each curve $\varphi|_{\partial D_j}$ is a Jordan curve and $$\sum_{j=1}^m \length(\varphi|_{\partial D_j}) \leq \length(c).$$
 \end{enumerate} 
\el
 
By a piecewise geodesic we mean the concatenation of finitely many geodesics, where we interpret constant curves as geodesics as well.  

\begin{proof}
We may assume that the points $s_0,\dots, s_N\in S^1$ are in cyclic order and we set $s_{N+1}=s_0$. Proceeding by induction on $i$ we find curves $\gamma _i\colon [s_0,s_i] \to X$ with the following properties.
The restriction $\gamma_i|_{[s_{k-1},s_k]}$ is a geodesic between $c(s_{k-1})$ and $c(s_k)$ for any $1\leq k\leq i$ and the image
$\gamma_i([s_0,s_i])$ is a finite topological graph $G_i$ in $X$ with geodesic edges. In the inductive step, we first choose an arbitrary geodesic $\eta$ from $c(s_i)$ to $c(s_{i+1})$.  After possibly modifying $\eta$ on finitely many intervals we may assume that $\eta$ intersects any edge of $G_i$ only at boundary points or connected subsets. Therefore the union of $G_i$ and the image of $\eta$ is a finite graph which has geodesic edges. We define $\gamma_{i+1}$ to be the concatenation of $\gamma_i$ and $\eta$.

For $i=N+1$ we obtain a piecewise geodesic $\gamma \colon S^1\to X$ with $\gamma(s_i)=c(s_i)$ for all $i$ and such that $\gamma(S^1)$ is a finite graph $G$ with geodesic edges. Moreover, by construction
$\length(\gamma) \leq \length(c)$. We may parametrize $\gamma$ to be Lipschitz continuous and define $\varphi|_{S^1}=\gamma$. Using  topological arguments in the finite graph $G$ we
easily extend 
$\varphi \colon S^1\to G$ to a Lipschitz continuous map $\varphi\colon A_m\to G$ such that (ii) holds true.  Since $\varphi$ has its image in $G$ the area of $\varphi$ is $0$.
\end{proof}

Let $X$ and $Y$ be metric spaces and $\varepsilon>0$. We say that $Y$ is an $\varepsilon$-thickening of $X$ if there exists an isometric embedding $\iota\colon X\to Y$ such that the Hausdorff distance between $\iota(X)$ and $Y$ is at most $\varepsilon$. The following lemma, which was proved in \cite{Wen08-sharp} and appeared in \cite[Lemma 3.3]{LWY16} in its present form, asserts the existence of $\varepsilon$-thickenings with good properties. 

\bl\label{lem:thickening-properties}
 Let $X$ be a length space.  There is a universal constant $M$ such that for every $\varepsilon>0$ there exists a complete length space $X_\varepsilon$ which is an $\varepsilon$-thickening of $X$ and has the following property. Let $\lambda>0$ and let $c_0\colon S^1\to X_\varepsilon$ be $\lambda$-Lipschitz.  If $\lambda \leq \frac{\varepsilon}{M}$, then $c_0$ is $M\lambda$-Lipschitz homotopic to a constant curve.  If $\lambda \geq \frac{\varepsilon}{M}$ then $c_0$ is Lipschitz homotopic to a curve $c_1\colon S^1\to X$ with $\length(c_1)\leq 2\length(c)$ via a homotopy of area at most $M \varepsilon \lambda$.  Furthermore, if $X$ is locally compact then $X_\varepsilon$ is locally compact.
\el

In particular, $X_\varepsilon$ is Lipschitz $1$-connected up to some scale.
We can now provide:

\begin{proof}[Proof of Theorem~\ref{thm:isop-ineq-for-geod-discs-not-optimal}]
 Let $k\geq 2$ and let $X_k$ be a $\frac{1}{k}$-thickening of $X$ as in Lemma~\ref{lem:thickening-properties}. Note that $X_k$ is locally compact and geodesic. Since $X_k$ converges in the Gromov-Hausdorff sense to $X$ as $k\to\infty$ it suffices, by \cite[Theorem 1.8]{LWY16}, to show that $X_k$ admits a quadratic isoperimetric inequality with isoperimetric constant only depending on $C$.  
 
Fix $k\geq 2$ and let $c\colon S^1\to X_k$ be a Lipschitz curve. We will show that $c$ extends to a Lipschitz map defined on $\overline{D}$ of area at most $C'\cdot \length(c)^2$, where $C'$ only depends on $C$.
By \cite[Lemma 3.6]{LWY16}, we may assume that $c$ has constant speed. By Lemma~\ref{lem:thickening-properties}, there exists a Lipschitz homotopy $\varphi_1$ of area  $$\Area(\varphi_1) \leq M\cdot \length(c)^2$$ from $c$ to a Lipschitz curve $c_1$ which is either constant or has image in $X$ and satisfies $\length(c_1)\leq 2\length(c)$. Here, $M$ denotes a suitable universal constant. If $c_1$ is constant then we are done, so we may assume $c_1$ to be non-constant and to have image in $X$. As above, we may assume that $c_1$ has constant speed. 

Fix $N\in\N$ sufficiently large, see below, and let $s_0, \dots, s_N\in S^1$ be equidistant points in cyclic order. Let $\varphi_2\colon A_m\to X$ be a Lipschitz map of zero area as in Lemma~\ref{lem:jordanize-curves}, when applied to the curve $c_1$. Set $c_2:= \varphi_2|_{S^1}$ and note that $$\length(c_2)\leq \length(c_1)\leq 2\length(c).$$ 
Set $s_{N+1}:= s_0$. If $N$ was chosen sufficiently large then there exist Lipschitz homotopies in $X_k$ of area at most $M'\cdot\frac{\length(c_1)^2}{(N+1)^2}$ from $c_1|_{[s_j, s_{j+1}]}$ to $c_2|_{[s_j, s_{j+1}]}$ for $j=0,\dots, N$, where $M'$ is a suitable universal constant. Using these homotopies we construct a Lipschitz homotopy $\varphi_2'$ from $c_1$ to $c_2$ satisying $$\Area(\varphi_2')\leq (N+1)\cdot M'\cdot\frac{\length(c_1)^2}{(N+1)^2} = M'\cdot \frac{\length(c_1)^2}{N+1}.$$

We consider the Lipschitz map $\varphi_2$ defined on $A_m = \overline{D}\setminus \cup_{j=1}^mD_j$. By construction, we have $$\sum_{j=1}^m\length(\varphi_2|_{\partial D_j})\leq \length(c_1)$$ and each $\varphi_2|_{\partial D_j}$ is a Jordan curve in $X$. Thus, by Proposition~\ref{prop:fill-Jordan-curve}, there exists a Lipschitz extension $\psi_j\colon \overline{D}_j\to X_k$ of $\varphi_2|_{\partial D_j}$  with $$\Area(\psi_j)\leq K\cdot\length(\varphi_2|_{\partial D_j})^2$$ for some constant $K$ depending only on $C$.

Finally, gluing the Lip\-schitz maps $\varphi_1$, $\varphi_2'$, $\varphi_2$, and $\psi_1,\dots, \psi_m$ we obtain a Lip\-schitz extension $\varphi$ of $c$ satisfying 
\begin{equation*}
 \begin{split}
  \Area(\varphi)&\leq \Area(\varphi_1) + \Area(\varphi_2') + \Area(\varphi_2) + \sum_{j=1}^m\Area(\psi_j)\\
  &\leq M\cdot \length(c)^2 + M'(N+1)^{-1}\cdot \length(c_1)^2 + K\cdot\sum_{j=1}^m \length(\varphi_2|_{\partial D_j})^2\\
  &\leq (M + 4M' + 4K)\cdot \length(c)^2.
 \end{split}
 \end{equation*}
This proves that $X$ admits a quadratic isoperimetric inequality with isoperimetric constant at most $C'= M + 4M' + 4K$.
\end{proof}

Theorem~\ref{thm:isop-ineq-for-geod-discs-not-optimal} has the following consequence:

\bc\label{cor:linloccon-isop-Sob}
 Let $X$ be a complete, geodesic metric space homeomorphic to $\overline{D}$, $S^2$, or $\R^2$. Suppose that $X$ is linearly locally connected and there exists $L>0$ such that $$\hm^2\left(B(x,r)\right)\leq L\cdot r^2$$ for every $x\in X$ and $r>0$. Then $X$ admits a quadratic isoperimetric inequality.
\ec

The corollary follows directly from Theorem~\ref{thm:isop-ineq-for-geod-discs-not-optimal} together with the lemma below.

\bl\label{lem:linloccon-settheoretic-isop}
 Let $X$ be a proper metric space homeomorphic to $\overline{D}$, $S^2$, or $\R^2$. Suppose that $X$ is linearly locally connected and there exists $L>0$ such that $$\hm^2(B(x, r))\leq L\cdot r^2$$ for every $x\in X$ and $r>0$. Then there exists $C>0$ such that every Jordan curve $\Gamma\subset X$ bounds a Jordan domain $\Omega\subset X$ satisfying $$\hm^2(\Omega)\leq C\cdot \length(\Gamma)^2.$$
\el

If $X$ is homeomorphic to $S^2$ or $\R^2$ then the lemma holds with $C=L\lambda^2$, where $\lambda$ is the linear local connectedness constant of $X$. If $X$ is homeomorphic to $\overline{D}$ then the constant $C$ which we obtain in our proof depends on $L$, $\lambda$, $\hm^2(X)$, and $\diam(\partial X)$.

\begin{proof}
We only give the proof in the case that $X$ is homeomorphic to $\overline{D}$, the argument for the other cases being similar but simpler. Let $\Omega\subset X$ be a Jordan domain and set $r:= \length(\partial \Omega)$. Let $\lambda\geq 1$ be the linear local connectedness constant. We distinguish two cases and first assume that $r<(2\lambda)^{-1}\cdot \diam(\partial X)$. Fix a point $x\in \partial \Omega$ and observe that there exists $x'\in \partial X$ with $d(x, x') > \lambda r$. Since $x'\not\in\Omega$ and $\partial \Omega\subset B(x,r)$ it follows that $x'\in X\setminus \overline{\Omega}$. We now show that $\Omega \subset B(x, \lambda r)$. We argue by contradiction and assume that there exists $x''\in \Omega\setminus B(x, \lambda r)$. Since $X$ is $\lambda$-linearly locally connected there exists a continuum $E\subset X\setminus B(x, r)$ connecting $x'$ and $x''$. However, $E$ must intersect $\partial \Omega$, which contradicts the fact that $\partial \Omega\subset B(x,r)$. Hence, $\Omega \subset B(x, \lambda r)$. Finally, we estimate $$\hm^2(\Omega)\leq \hm^2(B(x, \lambda r))\leq L\cdot (\lambda r)^2 = L\lambda^2\cdot\length(\partial \Omega)^2.$$ This concludes the proof of the first case. We now assume that $r\geq (2\lambda)^{-1}\cdot\diam(\partial X)$. Observe that $\hm^2(X)<\infty$ and hence $$\hm^2(\Omega) \leq \hm^2(X) \leq \frac{4\lambda^2\cdot \hm^2(X)}{\diam(\partial X)^2}\cdot \length(\partial \Omega)^2.$$ This concludes the proof. 
\end{proof}

We next establish Theorem~\ref{thm:isop-ineq-for-geod-discs}. The following proposition will be needed in its proof.

\bp\label{prop:almost-optima-filling}
 Let $X$ be as in Theorem~\ref{thm:isop-ineq-for-geod-discs} and let $c\colon S^1\to X$ be an injective Lipschitz curve. Then for every $\varepsilon>0$ there exists $u\in N^{1,2}(D, X)$ with $\trace(u) = c$ and such that $$\Area(u) \leq C\cdot\length(c)^2 + \varepsilon.$$
\ep

\begin{proof}
Let $c\colon S^1\to X$ be an injective Lipschitz curve. By \cite[Lemma 3.6]{LWY16}, we may assume that $c$ is parametrized proportional to arc-length.

Let $\Omega\subset X$ be the Jordan domain of smallest area and boundary $c$.
 Denote by $d_{\overline{\Omega}}$ the length metric on $\overline{\Omega}$. Set $Y:= (\overline{\Omega}, d)$ and $Z:= (\overline{\Omega}, d_{\overline{\Omega}})$, where $d$ is the metric from $X$. By Lemma~\ref{lem:change-metric}, the identity map $\iota\colon Z\to Y$ is a homeomorphism which preserves the lengths of curves and the Hausdorff $2$-measure of Borel subsets. Moreover, $Z$ is geodesic and $\iota$ is $1$-Lipschitz. Finally, $\iota$ is locally isometric on $Z\setminus \partial Z$. Hence, every Jordan domain $U\subset Z$ satisfies $$\hm^2_Z(U)\leq C\cdot \length(\partial U)^2.$$ It thus follows from Theorem~\ref{thm:isop-ineq-for-geod-discs-not-optimal} that $Z$ admits a quadratic isoperimetric inequality. In particular, $\Lambda(\partial Z, Z)$ is non-empty and, by Theorem~\ref{thm:existence-regularity-energy-min}, there exists $v\in \Lambda(\partial Z, Z)$ which minimizes the Reshetnyak energy $E_+^2$ among all maps in $\Lambda(\partial Z, Z)$. Moreover, $v$ has a unique representative which is continuous on $D$ and extends to a continuous map $v\colon \overline{D}\to Z$. Finally, $v$ satisfies Lusin's condition (N), see \cite[Theorem 4.4]{LW16-energy-area}. 
Theorem~\ref{thm:Plateau-solution-cell-like-intro} shows that $v$ is a cell-like map and thus monotone. 
 The area formula now implies that $$\Area(v) = \int_{v(D)} N(v,z) \, d\hm^2_Z(z),$$ where $N(v,z)$ denotes the number of points in the fiber $v^{-1}(z)$. Since $v$ is monotone and $\Area(v)<\infty$ it follows that $N(v,z)=1$ for almost every $z\in v(D)$ and hence
 $$ \Area(v)= \hm^2(Z)=\hm^2_Y(\overline{\Omega})\leq C\cdot \length(c)^2.$$
 
Let $\varepsilon>0$. By \cite[Lemma 4.8]{LW16-intrinsic}, we may connect $c$ and $\trace(v)$ by a Sobolev annulus $w$ of area at most $\varepsilon$. Gluing $v$ and $w$ we obtain a Sobolev map $u\in N^{1,2}(D, X)$ whose trace equals $c$ and with area at most $C\cdot \length(c)^2 + \varepsilon$. This completes the proof.
\end{proof}

We can finally provide:

\begin{proof}[Proof of Theorem~\ref{thm:isop-ineq-for-geod-discs}]
 Let $c\colon S^1\to X$ be a Lipschitz curve. 
 Fix $N\in\N$ sufficiently large, see below, and let $s_0,\dots, s_N\in S^1$ be equidistant points. Let $\varphi\colon A_m\to X$ be a Lipschitz map of zero area as in Lemma~\ref{lem:jordanize-curves}, when applied to the curve $c$. Set $c_1:= \varphi|_{S^1}$ and note that $\length(c_1)\leq \length(c)$. 
Exactly as in the proof of Theorem~\ref{thm:isop-ineq-for-geod-discs-not-optimal}, there exists a Lipschitz homotopy $\varphi'$ from $c$ to $c_1$ satisying $$\Area(\varphi')\leq M(N+1)^{-1}\cdot \length(c)^2$$ for some universal constant $M$, whenever $N$ was chosen large enough.

Let $\varepsilon>0$ and write $A_m$ as $A_m = \overline{D}\setminus \cup_{j=1}^mD_j$. By Proposition~\ref{prop:almost-optima-filling}, there exists for each $j$ a map $u_j\in N^{1,2}(D_j, X)$ with $\trace(u_j) = \varphi|_{\partial D_j}$ and $$\Area(u_j)\leq C\cdot \length(\varphi|_{\partial D_j})^2 + \frac{\varepsilon}{m}.$$ Gluing the maps $\varphi$, $\varphi'$, and $u_1,\dots, u_m$ yields a Sobolev map $u\in N^{1,2}(D,X)$ with $\trace(u) = c$ and such that $$\Area(u)\leq \Area(\varphi') +  \sum_{j=1}^m\Area(u_j)\leq \left[C + \varepsilon + M(N+1)^{-1}\right]\cdot \length(c)^2.$$
Choosing $\varepsilon>0$ arbitrarily small and $N$ arbitrarily large, we see that the isoperimetric constant of $X$ is at most $C$. This completes the proof.
\end{proof}

Combining Theorem~\ref{thm:isop-ineq-for-geod-discs}, Theorem~\ref{thm:existence-regularity-energy-min}, and Theorem~\ref{thm:Plateau-solution-cell-like-intro} we obtain the following almost parametrization result.

\bt\label{thm:almost-qc-param}
 Let $X$ be as in Theorem~\ref{thm:isop-ineq-for-geod-discs} and such that $X$ is homeomorphic to $\overline{D}$ and $\length(\partial X)<\infty$. Then $\Lambda(\partial X, X)$ contains an element $u$ of minimal energy $E_+^2(u)$. Every such $u$ is infinitesimally isotropic and has a continuous representative which is a uniform limit of homeomorphisms from $\overline{D}$ to $X$.
\et

The following example, which appeared in \cite[Example 11.3]{LW16-intrinsic}, illustrates that spaces as in the theorem need not be Ahlfors $2$-regular.

\begin{example}\label{ex:collapsed-disc}
 Let $T\subset D$ be a compact ball. Denote by $X$ the metric space obtained from $\overline{D}$ by identifying points in $T$, equipped with the quotient metric. Then $X$ is a geodesic metric space which is homeomorphic to $\overline{D}$ and satisfies \eqref{eq:geom-isop-intro}.
 \end{example}

Clearly, the space $X$ in the example is not Ahlfors $2$-regular. Moreover, $X$ is not reciprocal as defined in \cite{Raj14} and hence does not admit a quasiconformal parametrization in the sense of \cite{Raj14}.

\section{Proofs of parametrization results}\label{sec:proofs-main-thms}

The following result makes the statements in Theorem~\ref{thm:main-thm-qs-discs-summary} more precise and slightly more general.

\bt\label{thm:qs-subdomains-llc}
 Let $X$ be an Ahlfors $2$-regular, geodesic metric space homeomorphic to a $2$-dimensional manifold. Let $\Omega\subset X$ be a Jordan domain with $\length(\partial \Omega)<\infty$ and such that $\overline{\Omega}$ is linearly locally connected. Then there exists $u\in \Lambda(\partial \Omega, \overline{\Omega})$ which is continuous on $\overline{D}$ and satisfies $$E_+^2(u) = \inf\left\{E_+^2(v): v\in\Lambda(\partial \Omega, \overline{\Omega})\right\}.$$ Any such map is a quasisymmetric homeomorphism from $\overline{D}$ to $\overline{\Omega}$ and is uniquely determined up to a conformal diffeomorphism of $\overline{D}$.
\et

We first establish the following result:

\bt\label{thm:monotone-implies-homeo-quadgrowth}
 Let $X$ be a complete metric space and suppose there exists $L>0$ such that  for all $x\in X$ and $r>0$ we have $$\hm^2(B(x, r))\leq L\cdot r^2.$$ Let $u\colon \overline{D} \to X$ be continuous, monotone, and non-constant. If $u\in N^{1,2}(D, X)$ and $u$ is infinitesimally quasiconformal then $u$ is a homeomorphism onto its image.
\et

\begin{proof}
 It suffices to show that $u$ is injective. We argue by contradiction and suppose there exists $x\in u(\overline{D})$ such that $E:= u^{-1}(x)$ consists of more than one point. Fix $r>0$ such that $F_r:= \overline{D}\setminus u^{-1}(B(x,r))$ is not empty. Note that such $r$ exists since $u$ is assumed to be non-constant. Let $\Gamma:= \Gamma(E, F_r; \overline{D})$ denote the family of curves in $\overline{D}$ joining $E$ and $F_r$. 
 
We first show that the modulus of $\Gamma$ is bounded from above independently of $r$. Let $g_u$ be the minimal weak upper gradient of $u$ on $D$. It follows as in the proof of Proposition~\ref{prop:ae-qc-implies-weak-modulus-qc} that $g_u$ is a weak upper gradient of $u$ also on $\overline{D}$. We claim that the function $$\rho:= \frac{1}{r}\cdot g_u\cdot 1_{u^{-1}(B(x,r))}$$ is weakly admissible for $\Gamma$. Indeed, let $\gamma\colon [a,b]\to \overline{D}$ be a rectifiable curve in $\Gamma$ such that the upper gradient inequality \eqref{eq:upper-grad} holds with $g=g_u$ on $\gamma$ and on all its compact subcurves. This holds for almost every curve $\gamma\in\Gamma$ by \cite[Proposition 6.3.2]{HKST15}. Fix such $\gamma$. We may assume that $\gamma$ is parametrized by arc-length and satisfies $\gamma(a)\in E$. Let $t\leq b$ be the first point such that $d(u(\gamma(t)), x) = r$. Hence, by the upper gradient property, we have $$\int_\gamma \rho \geq  \int_a^t \rho(\gamma(s))\,ds \geq r^{-1}\cdot \int_a^t g_u(\gamma(s))\,ds \geq r^{-1}\cdot d(x, u(\gamma(t))) = 1.$$ This shows that $\rho$ is weakly admissible for $\Gamma$. 

Since $u$ is monotone and infinitesimally $Q$-quasiconformal for some $Q\geq 1$ we obtain, using the area formula as in the proof of Proposition~\ref{prop:ae-qc-implies-weak-modulus-qc}, that
\begin{equation*}
  \int_{D} \varrho^2(z)\,dz \leq \frac{Q}{r^2}\cdot \int_{u^{-1}(B(x,r))} \jac(\apmd u_z)\,dz\leq Q\cdot\frac{\hm^2(B(x,r))}{r^2}\leq QL.
 \end{equation*}
 This shows that $\MOD(\Gamma)$ is bounded from above independently of $r$. 
 
 Since $E$ is a non-degenerate continuum (i.e. consisting of more than one point), this bound on the modulus is easily seen to contradict the Loewner property of $\overline{D}$. Indeed, let $z_0\in \overline{D}\setminus E$ and let $z_1\in E$ be a point on $E$ nearest to $z_0$. For $s>0$ sufficiently small, let $G_s$ denote the set of points on the straight segment from $z_0$ to $z_1$ which are at least a distant $s$ away from $E$. Then $G_s$ is a non-degenerate continuum and, for $s>0$ sufficiently small, we have 
 $$\dist(E, G_s) \leq C\cdot s\cdot\min\{\diam E, \diam G_s\}$$ for some $C$ not depending on $s$. Fix $s$ as above and let $r>0$ be so small that $G_s\subset F_r$, where $F_r$ is as at the beginning of the proof. Then $\Gamma(E, F_r;\overline{D})$ contains the family $\Gamma(E, G_s; \overline{D})$ of curves in $\overline{D}$ connecting $E$ and $G_s$ and hence $$\MOD(\Gamma(E, G_s; \overline{D}))\leq \MOD(\Gamma(E, F_r;\overline{D})) \leq QL$$ for all $s>0$ sufficiently small. However, this is impossible since $$\MOD(\Gamma(E, G_s; \overline{D}))\to \infty$$ as $s\to 0^+$ by the $2$-Loewner property of $\overline{D}$, see e.g.~\cite[Theorem 8.23 and Example 8.24]{Hei01}. This completes the proof.
\end{proof}

We are ready for the proof of our main theorem concerning quasisymmetric parametrizations.

\begin{proof}[Proof of Theorem~\ref{thm:qs-subdomains-llc}]
Denote by $Y$ the set $\overline{\Omega}$ equipped with the metric from $X$ and by $Z$ the same set $\overline{\Omega}$ equipped with the length metric. By Lemma~\ref{lem:change-metric}, the identity map $\iota\colon Z\to Y$ is a homeomorphism which preserves the lengths of curves and the Hausdorff $2$-measure of Borel subsets. Moreover, $Z$ is geodesic and $\iota$ is $1$-Lipschitz. %Finally, if $v\colon D\to Z$ is continuous then $v\in N^{1,2}(D, Z)$ if and only if $\iota\circ v\in N^{1,2}(D, Y)$. In this case, $\apmd (\iota\circ v)_z = \apmd u_z$ for almost every $z\in D$ and thus $E_+^2(\iota\circ v) = E_+^2(v)$.

 By Lemma~\ref{lem:linloccon-settheoretic-isop}, there exists $C>0$ such that every Jordan domain $U\subset Y$ satisfies $$\hm^2(U)\leq C\cdot \length(\partial U)^2.$$ Hence, the same is true for Jordan domains in the space $Z$. Thus, Theorem~\ref{thm:isop-ineq-for-geod-discs-not-optimal} shows that $Z$ admits a quadratic isoperimetric inequality. It follows that also $Y$ admits a quadratic isoperimetric inequality. By Theorem~\ref{thm:existence-regularity-energy-min}, there exists $u\in \Lambda(\partial Y, Y)$ which minimizes the Reshetnyak energy $E_+^2$ among all maps in $\Lambda(\partial Y, Y)$. Moreover, any such $u$ is infinitesimally isotropic and has a unique representative which is continuous on $D$ and extends to a continuous map $u\colon \overline{D}\to Y$. 

We will now show that any $u$ with the properties above is a quasisymmetric homeomorphism. For this, consider the map $v:= \iota^{-1}\circ u$, which is continuous, in $N^{1,2}(D, Z)$, and satisfies  $\apmd v_z = \apmd u_z$ for almost every $z\in D$, see \cite[Corollary 3.2]{LW16-intrinsic}. In particular, $v\in \Lambda(\partial Z, Z)$ and $E_+^2(v) = E_+^2(u)$. It is now clear that $v$ is an energy minimizer in $\Lambda(\partial Z, Z)$ since for any $w\in \Lambda(\partial Z, Z)$ we have $$E_+^2(v) = E_+^2(u) \leq E_+^2(\iota\circ w) \leq E_+^2(w).$$
Thus, Theorem~\ref{thm:Plateau-solution-cell-like-intro} shows that $v$ is a uniform limit of homeomorphisms from $\overline{D}$ to $Z$ and thus monotone by Proposition~\ref{prop:equiv-cell-like-monotone}. Consequently, the map $u$ is monotone too. Since $u$ is infinitesimally isotropic and thus infinitesimally quasiconformal, Theorem~\ref{thm:monotone-implies-homeo-quadgrowth} shows that $u$ is a homeomorphism from $\overline{D}$ to $Y$. Proposition~\ref{prop:ae-qc-implies-weak-modulus-qc} and Theorem~\ref{thm:weak-modulus-qc-implies-qs} imply that $u$ is quasisymmetric.

We are left to show that any map as above is unique up to composition with a conformal diffeomorphism of $D$. Thus, let $u$ and $v$ be energy minimizers in $\Lambda(\partial Y, Y)$ which are continuous on $\overline{D}$. They are thus quasisymmetric homeomorphisms from $\overline{D}$ to $Y$ by the argument above. We will show that the map $\varphi\colon \overline{D}\to\overline{D}$ given by $\varphi:= v^{-1}\circ u$ is a  conformal diffeomorphism of $D$. First notice that, as the composition of two quasisymmetric homeomorphisms, the map $\varphi$ is itself quasisymmetric and, in particular, $\varphi$ and $\varphi^{-1}$ satisfy Lusin's condition (N), see \cite[Theorem 33.2]{Vai71}. It thus follows from the approximate metric differentiability of $u$ and $v$ that $$\apmd u_z = \apmd v_{\varphi(z)}\circ d_z\varphi$$ for almost every $z\in D$, see \cite[Proposition 4.3]{LW15-Plateau}. Here, $d_z\varphi$ denotes the classical differential of $\varphi$, which exists for almost every $z\in D$, see \cite[Corollary 32.2]{Vai71}. By the area formula and \cite[Corollary 8.15]{HKST01}, the approximate metric derivatives $\apmd u_z$ and $\apmd v_z$ are non-degenerate at almost every $z$. 

By Theorem~\ref{thm:existence-regularity-energy-min}, the maps $u$ and $v$ are infinitesimally isotropic, so for almost every $z$ the unique ellipses of largest area (John's ellipses) contained in the unit balls with respect to the norms $\apmd u_z$ and $\apmd v_{\varphi(z)}$ are given by Euclidean discs. Since the map $$d_z\varphi\colon (\R^2, \apmd u_z)\to (\R^2, \apmd v_{\varphi(z)})$$ is an isometry it must map John's ellipses to John's ellipses. This shows that $\varphi$ is $1$-quasiconformal and hence a conformal diffeomorphism. This proves that $u$ and $v$ indeed agree up to a conformal diffeomorphism.
\end{proof}

We can now prove a strenghtening of the Bonk--Kleiner theorem \cite{BK02}. Given a complete metric space $X$, denote by $N^{1,2}(S^2,X)$ the Newton-Sobolev space defined as in Section~\ref{sec:Sobolev-defs} with $\Omega$ replaced by $S^2$. The energy $E_+^2(u)$ of an element $u\in N^{1,2}(S^2, X)$ is defined analogously. Let $\Lambda(X)$ be the family of maps $u\in N^{1,2}(S^2, X)$ such that $u$ is a uniform limit of homeomorphisms from $S^2$ to $X$. 

\bt\label{thm:BK-canonical-version}
 Let $X$ be an Ahlfors $2$-regular,  linearly locally connected metric space homeomorphic to $S^2$. Then $\Lambda(X)$ is not empty and contains an element $u$ which satisfies $$E_+^2(u) = \inf\left\{E_+^2(v): v\in\Lambda(X)\right\}.$$ Any such map $u$ is a quasisymmetric homeomorphism from $S^2$ to $X$ and is uniquely determined up to a conformal diffeomorphism of $S^2$.
\et

We first show:

\bp\label{prop:Lambda-X-not-empty}
 If $X$ is as in Theorem~\ref{thm:BK-canonical-version} then $\Lambda(X)$ is not empty.
\ep

\begin{proof}
By \cite[Theorem B.6]{Sem96}, the space $X$ is quasi-convex. Thus, after changing the metric on $X$ in a biLip\-schitz way we may assume that $X$ is geodesic. Starting with a Jordan curve in $X$ and arguing exactly as in the proof of Lemma~\ref{lem:bilipap} we find a Jordan curve $\Gamma\subset X$ which is moreover a biLip\-schitz curve. Let $\Omega_{1,2}\subset X$ be the two Jordan domains enclosed by $\Gamma$. We claim that $\overline{\Omega}_1$ and $\overline{\Omega}_2$ are linearly locally connected. Indeed, one constructs continua in $\overline{\Omega}_{1,2}$ satisfying the linear local connectedness condition by using the continua in $X$ given by the linear local connectedness and replacing the part outside $\overline{\Omega}_{1,2}$ by a part of the biLipschitz curve $\partial \Omega_{1,2}$. 
By Theorem~\ref{thm:qs-subdomains-llc}, there exist quasisymmetric homeomorphisms $u_k\colon\overline{D}\to \overline{\Omega}_k$ for $k=1,2$. The map given by $\varphi:= u_2^{-1}\circ u_1|_{S^1}$ is a quasisymmetric homeomorphism of $S^1$ to itself and hence extends to a quasisymmetric homeomorphism $\overline{\varphi}$ of $\overline{D}$ by \cite{AB56}. Consequently, the quasisymmetric homeomorphism $\bar{u}_2\colon \overline{D}\to \overline{\Omega}_2$ defined by $\bar{u}_2:= u_2\circ\overline{\varphi}$ agrees with $u_1$ on the boundary $S^1$. Now, identify $S^2$ in a biLip\-schitz way with the space obtained by gluing two copies of $\overline{D}$ along their common boundary $S^1$. Let $\psi\colon S^2\to X$ be the homeomorphism which coincides with $\bar{u}_2$ on one copy of $\overline{D}$ and with $u_1$ on the other copy. Since $u_1$ and $\bar{u}_2$ are Sobolev maps it follows that $\psi$ is in $N^{1,2}(S^2, X)$ and hence $\Lambda(X)$ is not empty. This completes the proof. 
\end{proof}

Note that the map $\psi$ constructed in the proof satisfies the hypotheses of the quasisymmetric gluing theorem \cite[Theorem 3.1]{AKT05} and hence $\psi$ is a quasisymmetric homeomorphism. Since quasisymmetric homeomorphisms preserve the linear local connectednes we obtain, in particular, the Bonk--Kleiner theorem \cite{BK02} as a consequence.

\bc\label{cor:Bonk-Kleiner}
 Let $X$ be an Ahlfors $2$-regular metric space homeomorphic to $S^2$. Then $X$ is quasisymmetric to $S^2$ if and only if $X$ is linearly locally connected.
\ec

Now, we turn to the proof of Theorem~\ref{thm:BK-canonical-version}. We identify $S^2$ with the Riemann sphere $\hat{\mathbb C}$ and note that pre-compositions with conformal maps in domains of $\hat{\mathbb C}$ preserve Sobolev maps and the Reshetnyak energy, \cite{LW15-Plateau}.  This allows us to reduce all local questions and statements about elements in $N^{1,2} (S^2,X)$ to the case of Sobolev maps on bounded domains in $\mathbb C=\R ^2$. In particular, the Reshetnyak energy is lower semi-continuous for energy bounded sequences in  $N^{1,2} (S^2,X)$,  \cite{LW15-Plateau}, \cite{HKST15} and any map $u\in N^{1,2}(S^2,X)$ has an approximate metric derivative almost everywhere.  As in the case of discs we have:

\bt\label{thm:inf-isot-S2}
Let $u\ N^{1,2} (S^2,X)$ be such that  $E_+^2 (u\circ \phi) \geq E_+^2(u)$ for all biLipschitz homeomorphisms $\phi\colon S^2\to S^2$. Then 
$u$ is infinitesimally isotropic.
\et
Indeed, the proof of \cite[Theorem 1.2]{LW15-Plateau}, repeated and reformulated in  \cite[Lemma 3.2, Lemma 4.1]{LW16-energy-area}, applies to the present situation without change, since the ``critical'' biLipschitz homeomorphism $\phi$ of $\overline{D}$ constructed in \cite{LW15-Plateau} is fractional linear outside a small ball, hence extends to a biLipschitz homeomorphism of $\hat{\mathbb C}$ which is conformal outside a small ball.

\begin{proof}[Proof of Theorem~\ref{thm:BK-canonical-version}]
By Proposition~\ref{prop:Lambda-X-not-empty}, the family $\Lambda(X)$ is not empty. The existence of an energy minimizer in $\Lambda(X)$ is now proved as in the classical case when $X$ is smooth, see \cite[Section 3.1]{Jos91-book}. Indeed, fix distinct points $p_1, p_2, p_3\in S^2$ and distinct points $q_1, q_2, q_3\in X$. Let $(u_n)$ be an energy minimizing sequence in $\Lambda(X)$. After possibly composing with conformal diffeomorphisms of $S^2$ we may assume that each $u_n$ satisfies the $3$-point condition $u_n(p_i) = q_i$ for $i=1,2,3$. Applying the Courant-Lebesgue lemma, we deduce as in \cite[Section 3.1]{Jos91-book} that  the family $(u_n)$ is equi-continuous. Thus, after possibly passing to a subsequence, the sequence $(u_n)$ converges uniformly to a map $u\colon S^2\to X$. Then $u$ is in $N^{1,2}(S^2, X)$ and is a uniform limit of homeomorphisms from $S^2$ to $X$, hence $u\in \Lambda(X)$. By the lower semi-continuity of $E_+^2$ we thus see that $u$ is an energy minimizer in $\Lambda(X)$. This proves the existence of an energy minimizer in $\Lambda(X)$. 
 
Let $u$ be any energy minimizer in $\Lambda(X)$. By Theorem~\ref{thm:inf-isot-S2} the map $u$ is infinitesimally isotropic. Since $u$ is monotone it follows as in the proof of Theorem~\ref{thm:qs-subdomains-llc} that $u$ is a quasisymmetric homeomorphism and that $u$ is unique up to pre-composition with a conformal diffeomorphism.
\end{proof}

\section{Appendix}

The purpose of this section is to establish:

\begin{proof}[Proof of Theorem~\ref{thm:weak-modulus-qc-implies-qs}]
 Since quasisymmetric maps preserve linear local connectedness and the doubling property, one direction is clear. In order to prove the other direction suppose $X$ is doubling and linearly locally connected with some constant $\lambda\geq 1$. For $z\in \overline{D}$ and $r>0$ define $$L(z,r):= \max\{d(u(z), u(z')): z'\in \overline{D}\cap \bar{B}(z,2r)\},$$ where $\bar{B}(z,2r)$ denotes the closed ball. Since $u$ is uniformly continuous there exists $r>0$ such that 
\begin{equation}\label{eq:choice-of-r}
 L(z,r) < \frac{\diam X}{4}
\end{equation}
for all $z\in\overline{D}$.
 By \cite[Theorem 10.19]{Hei01} and \cite[Theorem 2.23]{TV80} it suffices to show that for every $z\in\overline{D}$ the restriction of $u$ to $\overline{D}\cap B(z, r)$ is a weak quasisymmetry. See \cite[Section 10]{Hei01} for the definition of a weak quasisymmetry.
 
 Fix $z\in\overline{D}$ and let $w,a,b\in \overline{D}\cap B(z,r)$ be such that $$|w-a|\leq |w-b|.$$ Set $s:= d(u(w), u(a))$ and let $M>0$ be such that $$s > M\cdot d(u(w), u(b)).$$ It suffices to show that $M$ must be bounded from above by a constant depending only on $Q$, $L$, and $\lambda$. If $M\leq 4\lambda^2$ then nothing needs to be proved, so we may assume that $M>4\lambda^2$. 
 
We claim that there exists $z'\in \overline{D} \setminus \bar{B}(z, 2r)$ with $u(z')\not\in \bar{B}(u(w), s/2)$. Indeed, due to \eqref{eq:choice-of-r}, there exists $y\in X$ with $d(y, u(z)) > 2L(z,r)$ and thus $z':= u^{-1}(y)$ satisfies $z'\not\in \bar{B}(z, 2r)$. Consequently, $$2L(z, r) < d(y, u(z))\leq d(y, u(w)) + d(u(w), u(z))\leq d(y, u(w)) + L(z, r),$$ from which it follows together with $$s\leq d(u(w), u(z)) + d(u(z), u(a)) \leq 2 L(z, r)$$ that $$d(y, u(w)) > L(z, r) \geq \frac{s}{2}.$$ This proves the claim.

Since $u(b)\in B(u(w), s/M)$ and $u(a), u(z') \not\in B(u(w), s/2)$ it follows from the linear local connectedness of $X$ that there exists a continuum $$E'\subset B(u(w), \lambda s/M)$$ connecting $u(w)$ and $u(b)$ and there exists a continuum $$F'\subset X \setminus B(u(w), s/2\lambda)$$ connecting $u(a)$ and $u(z')$.
 Let $\Gamma(E', F'; X)$ be the family of curves joining $E'$ to $F'$ in $X$. Since $\Gamma(E', F'; X)$ is a subset of the family of curves joining $\bar{B}(u(w), \lambda s/M)$ to $X\setminus B(u(w), s/2\lambda)$ in $X$ it follows from \cite[Lemma 7.18]{Hei01} that 
\begin{equation}\label{eq:ineq-balls-modulus} 
 \MOD(\Gamma(E', F'; X)) \leq L'\cdot \left(\log\frac{M}{2\lambda^2}\right)^{-1},
\end{equation} where $L'$ is a constant depending on $L$.

Set $E:= u^{-1}(E')$ and $F:= u^{-1}(F')$. Then $\Gamma(E', F'; X) = u\circ\Gamma(E, F; \overline{D})$ and hence 
 \begin{equation}\label{eq:ineq-mod-image}
  \MOD(\Gamma(E, F; \overline{D})) \leq Q\cdot \MOD(\Gamma(E', F'; X)).
 \end{equation}
We clearly have $$\dist(E, F)\leq 2\cdot \min\{\diam E, \diam F\}$$ since $\dist(E, F)\leq |w-a|$ and $$\min\{\diam E, \diam F\} \geq \min\{|w-b|, |a-z'|\}\geq \frac{1}{2}\cdot |w-a|.$$ It thus follows that 
\begin{equation}\label{eq:ineq-Loewner-proof}
 \MOD(\Gamma(E, F; \overline{D})) \geq \phi(2)>0,
\end{equation}
where $\phi$ is a Loewner function for $\overline{D}$. For a definition of the Loewner function and the Loewner property of $\overline{D}$ see for example \cite{Hei01}. Combining inequalities \eqref{eq:ineq-balls-modulus}, \eqref{eq:ineq-mod-image}, and \eqref{eq:ineq-Loewner-proof} we see that $$M\leq 2\lambda^2\cdot \exp\left(\frac{QL'}{\phi(2)}\right).$$This shows that the restriction of $u$ to $\overline{D}\cap B(z,r)$ is weakly $H$-quasisymmetric for some $H$ only depending on $\lambda, Q, L$.
\end{proof}

\def\cprime{$'$} \def\cprime{$'$} \def\cprime{$'$} \def\cprime{$'$}

\end{document}